\documentclass[final]{siamltex}

\usepackage{latexsym, enumerate}
\usepackage{eepic}
\usepackage{epic}
\usepackage{graphicx}
\usepackage{color}
\usepackage{ifpdf}
\usepackage{amssymb,amsmath,epsfig, algorithm,algorithmicx,multirow}
\usepackage{amsfonts, dsfont}
\usepackage{subfigure}
\usepackage{multirow}
\usepackage{makecell}
\usepackage{bm}

\newcommand{\bff}{\mathbf}
\newcommand{\bb}{\mathbb}
\newcommand{\ca}{\mathcal}

\title{Bayesian inference using intermediate distribution based on coarse multiscale model for time fractional diffusion equation
\thanks{L. Jiang acknowledges the support of Chinese NSF 11471107.
}}

\author{Lijian Jiang\thanks{Institute  of Mathematics, Hunan University, Changsha 410082, China. ({\tt  ljjiang@hnu.edu.cn}), Corresponding author.}
\and
Na Ou\thanks{College of Mathematics and Econometrics, Hunan University, Changsha 410082, China ({\tt oyoungla@163.com}).}
}
\begin{document}

\maketitle

\begin{abstract}
In the paper, we present  a strategy for accelerating posterior inference for unknown inputs in time fractional diffusion  models. In many inference problems, the posterior may be concentrated in a small portion of the entire prior support. It will be much more efficient if we  build and simulate  a surrogate only over the significant region of the posterior.
To this end, we  construct  a coarse model using Generalized Multiscale Finite Element Method (GMsFEM), and  solve a least-squares problem for the coarse model with a regularizing Levenberg-Marquart algorithm. An intermediate distribution is built based on the approximate sampling distribution.  For Bayesian inference,  we use GMsFEM and least-squares stochastic collocation method to obtain
a reduced coarse model based on the intermediate distribution. To increase the sampling speed of Markov chain Monte Carlo,  the DREAM$_\text{ZS}$ algorithm is used
to explore the surrogate posterior density, which is based on  the surrogate likelihood and the intermediate distribution.
The proposed method with lower gPC order  gives the approximate posterior as accurate as the  the surrogate model directly based  on the original prior.
 A few numerical examples for time fractional diffusion equations are carried out to demonstrate the performance of the proposed method with applications of the Bayesian inversion.
\end{abstract}



\pagestyle{myheadings}

\thispagestyle{plain}
\markboth{L. Jiang and Na Ou}{Bayesian inference for time fractional diffusion}

\section{Introduction}

Fractional-differential equations play a crucial role in  modeling  many physical phenomenon. They are broadly applied in damping laws, fluid mechanics, viscoelasticity, biology, physics, engineering and the modeling of earth quakes due to anomalous diffusion effects in constrained environments \cite{agrawal2004special, hilfer2000applications}. These fractional equations, which can be asymptotically derived from basic random walk models or generalized master equations, are physical-mathematical approach to anomalous diffusion. In practice,  the inputs or parameters in the  fractional models are often unknown and need to be identified based on some observation data and prior information.
The problem of identifying unknown inputs in mathematical models has been intensively  studied in the framework of  inverse problems and various numerical methods have been developed \cite{bondarenko2009numerical, miller2013coefficient, li2012numerical}.

Because this kind of inverse problems involve  sparse and indirect observations and the uncertainties from the prior information of  forward models, they are often ill-posed.
 In practical applications, the inevitable measurement error would render great  challenge in obtaining stable and accurate numerical solutions of the inverse problems. Iterative regularization methods \cite{doicu2010numerical} can be used to solve the ill-posed inverse problems. Another classical approach to regularize inverse problems is through the least-squares approach and penalty functional regularization \cite{engl2000regularization, tarantola2005inverse}. Both of them lead to point estimates of the parameters. However, what we are often interested in may not only point estimates but also the statistical distribution  of the parameters. This can be fulfilled by Bayesian inference. In this paper, we resort to Bayesian inference to identify unknown parameters in   fractional diffusion equations.

The Bayesian approach \cite{kaipio2006statistical, tarantola2005inverse, gamerman2006markov} incorporates uncertainties in noisy observations and prior information,  and can derive  the posterior probability density of the parameters, which enables us to quantify the uncertainty in the parameters. We can use the posterior conditional expectation or maximum a posterior (MAP) to characterize the parameters. However, the nonlinearity of the parameter-to-observation map and lack of analytical form of the forward model make the implementation of posterior expression very difficult. Alternatively, one can use Markov chain Monte Carlo (MCMC) \cite{liu2008monte, brooks2011Handbook, gamerman2006markov} techniques to construct Markov chains whose stationary distribution is the posterior density. The efficiency of many standard MCMC chains degrades with respect to the dimension of parameters \cite{roberts2002optimal, mattingly2012diffusion, pillai2012optimal}. The DREAM$_\text{ZS}$ \cite{vrugt2009accelerating, laloy2012high} algorithm
 employs simultaneous multiple Markov chains (generally 3-5) and uses the differential evolution algorithm to generate trial points for each chain. It
 has been particularly successful in finding appropriate trial parameters, which can improve the efficiency of MCMC chains.
 The DREAM$_\text{ZS}$ algorithm has successfully coped with high-dimensionality, nonlinearity, non-convexity, multimodality and numerous local optima problems \cite{vrugt2009accelerating}.

The MCMC approach entails repeated solutions of the forward model.  When the forward model is computationally intensive, such as multiscale models, a direct application of Bayesian inference would be computationally prohibitive. To significantly improve the simulation efficiency, reducing order or searching for surrogates of the forward models, or seeking more efficient sampling from the posterior \cite{efendiev2005efficient, martin2012stochastic} are necessary to accelerate the Bayesian inference. Numerical multiscale methods can efficiently and accurately solve multiscale  models in a coarse grid without resolving all scales in very fine grid.
Multiscale Finite Element Method (MsFEM) \cite{hou1997multiscale} is one of the  methods to efficiently simulate multiscale  models. The basic idea of MsFEM is to incorporate the small-scale information into multiscale basis functions and capture the impact of small-scale features on the coarse-scale through a variational formulation. One of the most important features for MsFEM
is that the multiscale basis functions can be computed overhead and used repeatedly for the model with different source terms and boundary conditions.
Many other multiscale methods \cite{allaire2005multiscale, e2003heterogeneous, hughes1998variational, jenny2003multiscale, jiang2012resourceful}  share the  similarities with MsFEM.
Generalized Multiscale Finite Element Method
(GMsFEM) \cite{efendiev2011multiscale, efendiev2013generalized} has been developed to solve multiscale models with complex multiscale structures. GMsFEM has some advantages over the standard MsFEM. For example, the coarse space in GMsFEM is more flexible and the convergence of GMsFEM is independent of the contrastness of the multiscales.

In Bayesian inversion, the model's output depends on the random parameters. We can use generalized polynomial chaos (gPC)-based stochastic Galerkin methods to propagate prior uncertainty through the forward model. As an alternative to the stochastic Galerkin approach, stochastic collocation requires only a few number of uncoupled deterministic simulations, with no reformulation of the governing equations
of the forward model. We assume that the model's output is a stochastic field and admits a gPC expansion. Then we choose a set of
collocation nodes and use least-squares methods to determine the coefficients of the gPC basis functions. We call the method as least-squares stochastic collocation method (LS-SCM). This method shares the same idea as probabilistic collocation method \cite{li2009stochastic}. LS-SCM inherits the merits from stochastic Galerkin methods and collocation methods. The recent work \cite{yan2015stochastic} employs a stochastic collocation algorithm using l1-minimization to construct stochastic sparse models with limited number of nodes, and their strategy has been applied to the Bayesian approach to handle nonlinear problems. The authors of \cite{jakeman2015enhancing} present a basis selection method that used with l1-minimization to adaptively determine the large coefficients of polynomial chaos expansions. Such sparse stochastic collocation methods may give a feasible approach to solve problems in high dimension random spaces.

The gPC surrogate is usually constructed based on a prior density \cite{jiang2016multiscale}. However, the posterior is concentrated in a small portion of the entire prior support in many inference problems. Thus, it may be much more efficient to build surrogates only over the important region of a posterior than the entire prior support. In this paper, we construct an intermediate distribution by solving a coarse reduced optimization problem, on which a surrogate is constructed. By the proposed  method,
we obtain the posterior density of unknowns, instead of getting point estimate or confidence intervals as the deterministic method do, the posterior density provides more quantities of interest.
 As a comparison to the two-stage MCMC method \cite{efendiev2005efficient}, which uses the coarse forward model to reject samplers with low probability via MH algorithm, the intermediate distribution can be easily obtained by a gradient-based iteration method, and it is not required to be close to the posterior or works as the proposal as in the two-stage MCMC method, i.e., the intermediate distribution is not a part of the Markov chain. Moreover,  the intermediate distribution is obtained from the sampling distribution and  considerably reduces uncertainty of the unknown parameters compared to the prior density, the gPC surrogate constructed based on the distribution leads to better approximate posterior density than the one derived by the original prior.  In the proposed approach, we  exclude the unimportant  region of the posterior by constructing the intermediate distribution, the acceptance rate of the Markov chains has been improved significantly.

Based on the intermediate distribution, we construct the reduced order model by combing GMsFEM with LS-SCM.  This  gives a representation for the model response. We solve an optimization problem based on the coarse reduced order model, where the measured data are used to inform the significant  region of the posterior, and we obtain the approximate sampling distribution. Then we  construct the intermediate distribution based on the sampling distribution, on which we use GMsFEM and LS-SCM to construct the surrogate model, and use the DREAM$_\text{ZS}$ algorithm to explore the surrogate posterior density. We have combined the deterministic and statistical method to solve inverse problems and obtain the statistical properties of the unknowns. We find that the surrogate based on the intermediate distribution and lower gPC order leads to the approximate posterior as accurate as the  the surrogate model  directly based on on the original prior. This can significantly alleviate the difficulty in simulating problems in high-dimensional  random spaces

The paper is organized as follows. We begin with formulating  a time-fractional diffusion equation and its inverse problems in Section 2 and introduce the framework of surrogate construction, especially the ideas of constructing the intermediate distribution.   Section 3 is devoted to reducing the dimensionality for the model parameter space  and the model state space. This includes  Karhunen-Lo\`{e}ve expansion, particularly  for the multiple random fields.  The model reduction method using GMsFEM and LS-SCM is presented in the section. In section 4, we present a few numerical examples to illustrate the performance of proposed method with applications in inverse time-fractional diffusion equation  problems. Some conclusions and comments are made finally.
%
%
%
\section{Statistical inversion based on intermediate distribution}

Anomalous diffusion has been well-known since Richardson's treatise on turbulent diffusion in 1926 \cite{richardson1926atmospheric}.
The role of anomalous diffusion has received attention within the literatures to describe many physical scenarios \cite{sun2013fractal, weiss2004anomalous}, most prominently within crowded systems, for example,  protein diffusion within cells and
diffusion through porous media. Fractional kinetic equations of the diffusion, diffusion advection and Fokker Planck type are presented as a useful approach for the description of transport dynamics in complex systems,
 which are governed by anomalous diffusion and non-exponential relaxation patterns.  The  typical  time-fractional diffusion equation is defined by
\begin{equation}
\label{frac-eq}
q(x)^c{D_t^\gamma} u-\text{div}(k(x)\nabla u)=f(x,t),\ \ x\in\Omega, t\in(0,T],
\end{equation}
subject to an appropriate boundary condition and initial condition.  Here $u=u(x,t)$ denotes state variable at space point $x$ and time $t$, and $\gamma \in (0, 1)$ is the fractional order of the derivative in time, where the Caputo fractional derivative of order $\gamma$ with respect to time is used and defined by
\begin{equation}
\label{cp-div}
^c{D_t^\gamma} u=\frac{1}{\Gamma(1-\gamma)}\int_0^{\text{t}} \frac{\partial u(x, s)}{\partial s} \frac{ds}{(t-s)^\gamma},
\end{equation}
where $\Gamma(\cdot)$ is the $\Gamma$ function. Here $k(x)$ is a spatial varied diffusion coefficient, $q(x)$ is the specific field and  the term $f(x, t)$ is a source (or sink) term. To simplify the function notations, we will suppress the variables $x$ and $t$ in functions when no ambiguity occurs. For practical models,  the model inputs such as $q(x)$, $k(x)$ and $\gamma$  may be not known, and they need to be estimated by some observations or measurements. In the paper, we focus on the inverse problems for the fractional diffusion model (\ref{frac-eq}).

\subsection{Bayesian inference}

In the paper, we use Bayesian inference to estimate unknown parameters  by some given noisy measurements of the model response at various sensors. We consider the case of additive noise $\varepsilon$ with probability density function $\pi_N(\varepsilon)$, then the measurement data can be expressed by
\begin{equation}
\label{stac-mod}
D = \bff{G(m)}+\varepsilon,
\end{equation}
where $\bff m$ is a vector of model parameters or inputs and $\bff{G(m)}\in\bb{R}^{n_d}$ is the observed model response at measurement sensors. Here  $n_d$ is the dimension of observations. We assume that $\varepsilon$ is independent of $\bff m$, then the conditional probability density for the measurement data $D$ given the unknown $\bff m$, i.e., the likelihood function is given by
\begin{equation}
\label{likeli_noise}
\pi(d|\bff m)= \pi_N \big(d-\bff{G(m)}\big).
\end{equation}
We use Bayesian inference to solve the inverse problem. This approach results in measure of uncertainty but not only a single-point estimate of the parameter. This is an advantage of Bayesian method over the standard regularization method. In the Bayesian setting, both $\bff m$ and $D$ are random variables. Then the posterior probability density for $\bff m$ can be derived by  Bayesian rule,
\begin{equation}
\label{Bayes}
\pi(\bff m|d)\propto \pi(d|\bff m)\pi(\bff m),
\end{equation}
where $\pi(\bff m)$ is the prior distribution with available prior information before the data is observed. The  data enters the Bayesian formulation through the likelihood function $\pi(d|\bff m)$. For the convenience of notation, we will use $\pi^d(\bff m)$ to denote the posterior density $\pi(\bff m|d)$ and $\bff{L}(\bff m)$ to denote the likelihood function $\pi(d|\bff m)$. Then (\ref{Bayes}) can be rewritten as
\begin{equation}
\label{like-L}
\pi^d(\bff m)\propto \bff{L}(\bff m)\pi(\bff m).
\end{equation}
Furthermore, if the prior density is conditional to unknown parameter $\theta$, i.e., $\pi(\bff m|\theta)$, the parameter $\theta$ is also a part of the inference problem in the Bayesian framework. In other words, these hyperparameters may be endowed with priors and estimated from data
\[
\pi(\bff m, \theta|d)\propto \bff{L}(\bff m)\pi(\bff m|\theta)\pi(\theta).
\]
The vector $\varepsilon$ is assumed to be independent and identically distributed (i.i.d.) Gaussian random vector with mean zero and standard deviation $\sigma_0$,
\[
\varepsilon\sim N(0,\sigma_0^2\mathbb{I}),
\]
where $\mathbb{I}$ is the identity matrix of size $n_d\times n_d$. Then the likelihood $\bff{L}(\bff m)$ defined in $(\ref{likeli_noise})$ is given by
\begin{equation}
\label{LH}
\bff{L}(\bff m) = (2\pi\sigma_0^2)^{-\frac{n_d}{2}}\exp\bigg(-\frac{\|d-\bff{G(m)}\|_2^2}{2\sigma_0^2}\bigg),
\end{equation}
where $\| \cdot \|_2$ refers to the Euclidean norm. We note that it is not necessary to compute the normalized term in (\ref{like-L}) under most circumstances. As the posterior distribution of $\bff m$ can be inferred, we can extract the posterior mean or the maximum a posteriori (MAP) of the unknowns
\[
\bff m_\text{MAP}=\arg\max \pi^d(\bff m).
\]
The MAP estimate is equivalent to the solution of a regularization minimization problem for some specific priors. However, the analytical expression of the posterior distribution is generally unavailable and the high dimension integration involved in posterior expectation is a great challenge. Markov chain Monte Carlo (MCMC) methods are a class of algorithms for sampling from a probability distribution based on constructing a Markov chain that has the desired distribution as its equilibrium distribution, and we can use the MCMC methods to explore the posterior state space of the unknowns.

We can calculate the credible intervals for inferred parameters and the model response based on posterior samplers. The $(1-\alpha_0)\times 100\%$ credible interval for a random variable is defined as
\[
\text{Prob}(\mu_L<\bff m<\mu_R)=1-\alpha_0,
\]
and the $(1-\alpha_0)\times 100\%$ prediction interval for a random response $\bff{u^\delta}(x,t)$ is referred to the pair of statistics $[\Upsilon_L, \Upsilon_R]$ constructed from a random observation $D$ such that
\[
\text{Prob}(\Upsilon_L<\bff{u^\delta}(x,t)<\Upsilon_R)=1-\alpha_0,
\]
where $\bff{u^\delta}(x,t)$ is a new observation at the point $(x,t)$, which  is independent of the data used to construct $\Upsilon_L$ and $\Upsilon_R$. Once the corresponding samplers obtained, the credible interval and prediction interval are easily constructed. For example,  for $\alpha_0=0.05$, we put the samplers in ascending order and determine the location of the 0.025 and 0.975 values, since we have the Monte Carlo integration in calculating the $p$ quantile $\mu_0$ for a random variable as
\[
p=\int_{-\infty}^{\mu_0} \pi^d(\bff m)d\bff m
=\frac{1}{M}\sum_{k=1}^M \bff{I}(\bff m^k<\mu_0)
=\frac{M_{\mu_0}}{M},
\]
where the sequence $\{\bff m^k\}_{k=1}^M$ are the samplers taken from the density $\pi^d(\bff m)$, and $\bff{I}( \cdot )$ is an indicator function. We can obtain the location of $p$ by $M_{\mu_0}=pM$, i.e., the $M_{\mu_0}$th sampler in the ascended samplers is the $p$ quantile $\mu_0$.

\subsection{Reduced order model  and intermediate distribution}

Let $\ca{M}$ be parameter space, and $\bff{m}\in \ca{M}$, we can directly apply the Bayesian formulation described in the previous section and explore the posterior density of $\bff{m}$ with MCMC. When using finite element method (FEM) to solve the forward model, the equation (\ref{frac-eq}) leads to a system of algebraic equations
\[
\bff{K(m)u}=\bff{f},\ \ \bff{G}(\bff{m})=\bff{Cu},
\]
where $\bff{K(m)}\in \bb{R}^{N_h\times N_h}$ is the nonlinear forward operator depending on the parameter $\bff{m}\in \bb{R}^n$, $\bff{u}\in \bb{R}^{N_h}$ is the solution or state, $\bff{f}\in \bb{R}^{N_h}$ is the source and $N_h$ is the degree of freedoms in the FEM, and $\bff{C}\in \bb{R}^{n_d\times N_h}$ is the observation operator, $\bff{G(m)} \in \bb{R}^{n_d}$ is the observation vector.

The dimension of the random vector $\bff{m}$ may  be high-dimensional because it  depends  on the  discretization of the physical domain. High dimensionality will make  MCMC exploration of posterior more challenging. Besides, we note that the dimension number of the states scales with the degree of freedoms of the FEM, but the dimension of the observable outputs $n_d \ll N_h^2$ in practice. Thus, it is necessary to construct a reduced model to accelerate the forward model computation. We assume that the
parameter $\bff{m}$ and state $\bff{u}$ can be adequately approximated in the span of parameter and state bases $\bff{E}\in \bb{R}^{n\times n_z}$ and $R\in \bb{R}^{N_h\times M_v}$, respectively.  This leads to  a reduced  algebraic system corresponding to (\ref{frac-eq}) as follows
\[
\bff{K_r}(z)\bff{u_r}=\bff{f_r},\ \ \bff{G_r}(z)=\bff{C_ru_r},
\]
where
\begin{equation}
\label{reduce}
\bff{K_r}(z)=R^T\bff{K}(\bff{E}z)R, \ \ \bff{f_r}=R^T\bff{f}, \ \ \bff{C_r}=\bff{C}R.
\end{equation}
Here $\bff{K_r}(z)\in \bb{R}^{M_v\times M_v}$ denotes  the reduced forward model depending on the reduced parameter
$z\in \bb{R}^{n_z}$, $\bff{u_r}\in \bb{R}^{M_v}$ is the reduced state, $\bff{f_r}\in \bb{R}^{M_v}$ is the projected source, $\bff{C_r}\in \bb{R}^{n_d\times M_v}$ is the reduced model observation operator, and $\bff{G_r}(z)\in  \bb{R}^{n_d}$ are the reduced model outputs. The reduction in parameter space is enforced directly by assuming $\bff{m}= \bff{E}z$. We note $M_v \ll N_h$ and $n_z\ll n$. Thus, both  the states and the parameters reside in lower dimensional spaces.

In additional to  reducing the parameter space and state space, we also use stochastic response surface methods (e.g., generalized polynomial chaos)  to construct surrogate, which can be used to accelerate evaluations of the posterior density. Generally, we can construct the surrogate based on the prior density $\pi(z)$ using the reduced model and stochastic response surface methods.  However, constructing a sufficiently accurate surrogate model over the support of the prior distribution may not be possible in many practical problems. As the posterior density reflects some information gain relative to the prior distribution, this motivates  us to seek  an intermediate distribution from a sampling distribution. To this end, we solve the  nonlinear optimization problem
for the derived  reduced order model
\begin{equation}
\label{coarseopt}
\mathcal{F}_r(z)=\min\{\frac{1}{2}\|d-\bff{G_r}(z)\|_2^2\}.
\end{equation}
From a Frequentist perspective, given the observed data, the solution to problem (\ref{coarseopt}) is one realization of the ordinary least square (OLS) estimator $\hat z$, where
\[
\hat z=\arg\min \frac{1}{2}\|D-\bff{G_r}(z)\|^2_2.
\]
Here $D\in \bb{R}^{n_d}$ is a random vector as noted in equation (\ref{stac-mod}), $d$ is a realization of $D$. Let  $z_0$ denote the true but unknown value of the parameter set that generated the observation $d$. Then we have the following conclusion \cite{smith2013uncertainty}: when the measurement errors are i.i.d. and $\varepsilon\sim N(0, \sigma_0^2)$, or $n_d$ is sufficiently large, we can specify a sampling distribution for $\hat z$ and it can be directly or asymptotically established that
\[
\hat z \sim N\bigg(z_0, \bff V\bigg),
\]
where the covariance matrix is given by
\[
\bff V=\sigma_0^2[\bff{J}^T(z_0)\bff{J}(z_0)]^{-1},
\]
and $\bff{J}(z)$ denotes the $n_d\times n_z$ sensitivity matrix whose elements are
\[
\bff{J}_{ik}(z)=\frac{\partial [\bff{G_r}]_i(z)}{\partial z_k}.
\]
The asymptotical distribution is the so-called sampling distribution. When we use the reduced model to approximate the forward model as shown in the optimization problem (\ref{coarseopt}), the model reduced error is unknown.
 This  leads the error between measured data and the reduced model unknown, the covariance matrix $\bff V$ can be approximated by
\[
\bff V\approx \sigma_\text{OLS}^2[\bff{J}^T(z_\text{OLS})\bff{J}(z_\text{OLS})]^{-1},
\]
where
\[
\sigma_\text{OLS}^2 = \frac{1}{n_d-n_z}\bff R^T\bff R, \quad
\bff R= d-\bff{G_r}(z_\text{OLS}).
\]
Furthermore, if we use $\delta_k$ to denote the $k-$th diagonal element of $[\bff{J}^T(z_\text{OLS})\bff{J}(z_\text{OLS})]^{-1}$ and $z_{0_k}$  the $k-$th element of the true parameter $z_0$, then the $(1-\alpha_0)\times 100\%$ confidence interval for $z_{0_k}$ is
\[
[z_\text{OLS$_k$} -t_{n_d-n_z, 1-\alpha_0/2}\sigma_\text{OLS}\sqrt \delta_k, z_\text{OLS$_k$} +t_{n_d-n_z, 1-\alpha_0/2}\sigma_\text{OLS}\sqrt \delta_k].
\]
As $n_d\gg n_z$ in practice, the t-distribution with $n_d-n_z$ degrees of freedom is approximated by the normal distribution. Following the $2\sigma$ rule of normal distribution, we have the assertion: for each $z_{0_k}$, it has high probability to be included in the interval
\[
[z_\text{OLS$_k$} -2\sigma_\text{OLS}\sqrt \delta_k, z_\text{OLS$_k$} +2\sigma_\text{OLS}\sqrt \delta_k],
\]
where $z_\text{OLS}$ is a realization of $\hat z$. We have gotten an interval estimation for $z_0$ in the probability sense, and  assume it uniformly distributed on the interval, which is used  as the new prior in Bayesian inference. To distinguish it from the original prior,  we call the new prior as  the intermediate distribution. The diagonal element of $\bff{V}$ is often very small, which implies the support of the intermediate distribution narrower than the original prior's. We will construct surrogate based on the intermediate distribution.

\subsection{DREAM$_\text{ZS}$ sampling method}

After obtaining the approximate posterior for the reduced parameter $z$, which is cooperated by the surrogate likelihood and the intermediate distribution, we use DREAM$_\text{ZS}$ algorithm to explore the posterior density and obtain samplers. DREAM$_\text{ZS}$ has been widely used  to  increase the speed of the MCMC process. It employs multiple Markov chains (generally 3-5) simultaneously and uses a differential evolution algorithm to generate trial points for each chain.

The DREAM$_\text{ZS}$ algorithm stems from DREAM but creates an {\em archive} of potential solutions.  Thus  we can  focus on the introduce to  the DREAM algorithm. For $n_z$-dimension parameters, $N_C$ chains $\{z^i\}$ ( $i=1,\cdots ,N_C$) simultaneously run in parallel and the present population is stored in an $N_C \times n_z$ matrix $X$. The candidates are randomly generated according to
\begin{equation}
\label{proposal}
\bff{q}^i=z^i+(\bff{1}_{n_z}+\bff{e})g'(\beta,n'_z)[\sum_{j=1}^\beta z^{r_1(j)}-\sum_{j=1}^\beta z^{r_2(n)}]+\bff{\epsilon},
\end{equation}
where $\beta$ signifies the number of pairs used to generate the proposal, and $r_1(j),  r_2(n)\in \{1,\cdots,N_C\}$, $r_1(j)\neq r_2(n)\neq i$ for $j, n=1,\cdots,\beta$. The value of $\bff{e}$ and $\bff{\epsilon}$ are drawn from $U(-h_1,h_1)$ and $N(0,h_2)$ with $|h_1|<1$ respectively, the $h_2$ is small compared to the width of the target distribution. The value of jump-size $g$ depends on $\beta$ and $n'_z$, i.e., $g(\beta,n'_z)=2.38/\sqrt{2\beta n'_z}$ and we set $g(\beta, n'_z)=1$ at every 5th generation.  The $g'$ is described in table \ref{dreamtab}.

Subspace sampling is implemented in DREAM by only updating randomly selected dimensions of $z^i$ each time a proposal is generated.
If $A_s$ is a subset of $n'_z$-dimensions of the original parameter space, $\bb{R}^{n'_z}\subseteq \bb{R}^{n_z}$, then a jump in the $i-$th chain at iteration $k+1$ is calculated using differential evolution crossover operation. The crossover operation is used to decide whether to update the $j-$th dimension of the parameter and the selection process follows the rule (for chain $i$ and $j-$th parameter)
\[
\bff{q}^i_j=\begin{cases}
z^i_j &\text{if}\ $U$\leq $1-CR$\\
\bff{q}^i_j &\text{otherwise},
\end{cases}
\]
where $\text{CR}$ denotes the crossover probability, and $U\in[0,1]$ is a draw from a uniform distribution. If the user-defined parameter is $n_\text{CR}$, then $\text{CR}$ is taken from the multinomial distribution. The pseudo-code of DREAM$_\text{ZS}$ algorithm is presented in Algorithm \ref{dreamtab}.
\begin{algorithm}
\caption{The DREAM$_\text{ZS}$ algorithm to sample $\tilde \pi^d(z)$}
      Initialise the {\em archive} with length $N_0$ and draw an initial population $X$ using the prior distribution;\\
      \textbf{while $j<N_E$}\\
      \textbf{for $i=1:N_C$} \\
      Compute the density $\tilde \pi^d(X^i)$;\\
      Generate a candidate point, $\bff{q}^i$ as (\ref{proposal}), where the differential evolution pairs are taken from the {\em archive} ; \\
      Replace each element $(j=1,\cdots,n_z)$ of the proposal $\bff{q}_j^i$ with $z_j^i$ using a binomial scheme with probability 1-CR;\\
      Compute $\alpha(X^i, \bff{q}^i)=\min\{1, \frac{\tilde\pi^d(\bff{q}^i)}{\tilde\pi^d(X^i)}\}$;\\
      \textbf{if rand(1)$<\alpha$} then\\
      Set $X^i=\bff{q}^i$;\\
      \textbf{end if}\\
      \textbf{end for}\\
      Store the samplers: $S_z=[S_z;X]$;\\
      Update the {\em archive}: {\em archive}$=[X;archive]$;\\
      Update the length of the {\em archive} to be $N_0+N_C$;\\
      \textbf{end while}\\
    \label{dreamtab}
\end{algorithm}

\section{Construction of surrogate model}
\label{method}

As we can see, obtaining samplers by  DREAM$_\text{ZS}$ algorithm requires us to run 3-5 chains simultaneously, which means we need to call a large numbers of  deterministic forward solvers. This may be computationally expensive and inefficient. We appeal to the truncated KLE and GMsFEM to reduce the dimension of unknown parameter space and state space  in equation (\ref{reduce}).
We assume that
\[
\bff{E}=[\sqrt{\lambda_1}e_1(x),\cdots, \sqrt{\lambda_{n_z}}e_{n_z}(x)],
\]
is the weighted eigenfunctions in KLE, and $e_{i}(x)\in \bb{R}^n$, $i=1,\cdots, n_z$, $z \in \bb{R}^{n_z}$ is the KLE modes, and the matrix $R$ is the multiscale basis. We note that GMsFEM provides us much flexibility in choosing multiscale basis functions, so we can set a small  number of selected basis functions on each coarse neighborhood  while maintaining the accuracy to some degree. Hence, we set the number of multiscale basis as $M_c$ in each coarse block to construct the reduced observation vector $\bff{G_r}(z)$, and use the regularizing Levenberg-Marquart algorithm to solve the ill-posed nonlinear problem (\ref{coarseopt}).
The computational cost of solving PDE involved in optimization iterations can be reduced significantly
 due to the scale reduction of the forward model. The iteration associated with regularizing Levenberg-Marquart algorithm is given by
\[
z^{k+1}=z^k+J^\dag \bigg(d-\bff{G_r}(z^k)\bigg),
\]
where
\[
J^\dag=(\bff{J}^T\bff{J}+\alpha \bb{I})^{-1}\bff{J}^T.
\]
Here the $\alpha$ is the regularization parameter and $\bb{I}$ is the identity matrix of size $n_z\times n_z$, $\bff{J}$ is the sensitivity matrix, whose transpose is defined by
\[
\bff{J}^T(z^k)=\bigg [\nabla \bff{G_r}_1(z^k), \nabla \bff{G_r}_2(z^k), \cdots, \nabla \bff{G_r}_{n_d}(z^k)\bigg ].
\]
The gradient of $\bff{G_r}_i(z^k)$,  $i=1,\cdots ,n_d$ can be approximated by a difference method, e.g.,  each column of $\bff{J}$ can be computed  by
\[
\bff{J}(:,j)=\frac{\bff{G_r}(z^k+h\eta_j)-\bff{G_r}(z^k)}{h}
\]
for $j=1,2,\cdots,n_z$, where $h$ is the stepsize and $\eta_j$ is a vector of zeros with one in the $j-$th component.

In reality, the accurate solution to problem (\ref{coarseopt}) is unaccessible, both the discretization of the forward model and the iterative methods may cause errors.
This would cause  some difference from the accurate interval estimation of the true parameter $z_0$. Thus we need some adjustment in determining the radius of the interval.
We recall that the interval estimation of $z_0$ is
\[
[z_\text{OLS$_k$}-2\sigma_\text{OLS}\sqrt\delta_k, z_\text{OLS$_k$}+2\sigma_\text{OLS}\sqrt\delta_k].
\]
Because the diagonal element of $\bff{V}$ is often very small,  the support of the intermediate distribution would be small and
we need to do some extension  to avoid the concentration at a very narrow region. We denote the modified  intermediate distribution of $z_{0_k}$ by $U(z_\text{OLS$_k$}-\nu, z_\text{OLS$_k$}+\nu)$, where the radius $\nu$ is chosen following the rule: if $2\sigma_\text{OLS}\sqrt\delta_k\leq 1$, we set $\nu=1$, else $\nu$ is set as $2\sigma_\text{OLS}\sqrt\delta_k$. The adjustment may not be optimal, we choose 1 as the criterion here just because we want to use Legend polynomials to construct surrogate in computation.  We note that the support of the modified intermediate distribution is much narrower than $\bb{R}$.

In summary, we use the regularizing Levenberg-Marquart algorithm to solve the problem (\ref{coarseopt}), and construct an intermediate distribution, on which to build the surrogate using LS-SCM and  GMsFEM. Then we  obtain an approximate posterior by integrating  the surrogate likelihood and the intermediate distribution. Furthermore, we explore the approximate  posterior by MCMC method and characterize the statistical properties of unknowns by samplers.  The outline of surrogate construction is described in Table \ref{tab2}.  We can use the Kullback-Leibler divergence to measure the difference between the approximate posterior and the full posterior.
\begin{table}[htbp]
  \centering
  \caption{The outline of surrogate construction}\label{tab2}
  \begin{tabular}{l}
     \hline
    $\bullet$ Calculate the GMsFEM matrix $R$; \rule{0pt}{0.5cm}\\
    $\bullet$ Solve the reduced optimization problem and obtain the solution $z_\text{OLS}$ \rule{0pt}{0.5cm}\\
    \quad \quad \quad \quad \quad $\mathcal{F}_r(z)=\min\{\frac{1}{2}\|d-\bff{G_r}(z)\|_2^2\}$;\rule{0pt}{0.8cm}\\
    $\bullet$ Calculate the covariance matrix $\bff{V}$ and construct the initial intermediate distribution as \\
    \quad \quad \quad \quad \quad $U(z_\text{OLS$_k$}-\nu, z_\text{OLS$_k$}+\nu)$;\rule{0pt}{0.8cm}\\
    $\bullet$ Do adjustment to the intermediate distribution:\rule{0pt}{0.5cm}\\
    \quad \quad \quad\textbf{for} $k=1:n_z$\rule{0pt}{0.5cm}\\
    \quad \quad \quad \textbf{if} $2\sigma_\text{OLS}\sqrt\delta_k\leq 1$\rule{0pt}{0.5cm}\\
    \quad \quad \quad \quad $\nu=1$;\rule{0pt}{0.5cm}\\
    \quad \quad \quad \textbf{else} $2\sigma_\text{OLS}\sqrt\delta_k>1$\rule{0pt}{0.5cm} \\
    \quad \quad \quad \quad $\nu=2\sigma_\text{OLS}\sqrt\delta_k$;\rule{0pt}{0.5cm}\\
    \quad \quad \quad \textbf{end if}\rule{0pt}{0.5cm}\\
    \quad \quad \quad \textbf{end for}\rule{0pt}{0.5cm}\\
    $\bullet$ Construct surrogate based on the intermediate distribution using LS-SCM and finer \rule{0pt}{0.5cm}\\
    \ \ GMsFEM;\rule{0pt}{0.5cm}\\
    $\bullet$ Explore the approximated posterior using DREAM$_\text{ZS}$ algorithm.\rule{0pt}{0.5cm}\\
     \hline
   \end{tabular}
\end{table}

\subsection{Karhunen-Lo\`{e}ve expansion}

Let $Y_1(x,\omega)$ and $Y_2(x,\omega)$ be two Gaussian random fields with second-order moments. We first parameterize the random  fields $Y_i(x,\omega)$ ($i=1,2$) by applying  Karhunen-Lo\`{e}ve expansion (KLE)  with given covariance functions and  truncate the KLE to approximately represent the random  fields.

For a general presentation, we consider multiple random input fields $Y_i(x,\omega)$ $(i=1,\cdots,m)$ defined on the same physical domain $\Omega$, which can be normalized as $Y'_j(x,\omega)=\big(Y_j(x,\omega)-\mathbb{E}[Y_j(x,\omega)]\big)/\sigma_{Y_j}$, where $\sigma_{Y_j}$ is the standard deviation of $Y_j$ and $\bb{E}[\cdot]$ is the expectation operator. For a single random field $ Y'_j(x,\omega)$, it can be expressed as
\[
Y'_j(x,\omega)=\sum_{i=1}^\infty \sqrt{\lambda_i^j}e_i^j(x)Z^j_i(\omega),
\]
where $\{Z^j_i(\omega)\}_{i=1}^\infty$ are uncorrelated random variables, $\{\lambda^j_i\}_{i=1}^\infty$ and $\{e^j_i(x)\}_{i=1}^\infty$ are the eigenpairs defined by
\[
\int_\Omega \text{cov}[Y'_j](x_1,x_2)e^j(x_2)=\lambda^je^j(x_1),
\]
where $\text{cov}[Y'_j]$ is the covariance  function of $Y'_j(x,\omega)$. In particular, if $Y'_j(x,\omega)$ is a Gaussian field, then $\{Z^j_i(\omega)\}_{i=1}^\infty$ are standard identically independent Gaussian random variables. We sort the eigenvalues in ascending order, i.e.,  $\lambda_1^j\geq\lambda_2^j\cdots$, and their corresponding eigenfunctions are also sorted accordingly. The $ Y'_j(x,\omega)$ is approximated by the truncated KLE expansion with the first $N_j$ terms, which  depends on the spectrum energy measured by $\sum_{i=1}^{N_j}\lambda_i^j/|\Omega|$, where $|\Omega|$ is the area of the spatial domain $\Omega$.

For multiple random input fields, the correlation between them can be described by the symmetric and positive definite covariance matrix with the $(i,j)$th element $\rho_{i,j}$,
\[
\rho_{i,j}=\bb{E}[Y'_i(x,\omega)Y'_j(x,\omega)],
\]
We decompose the covariance matrix into $LL^T$ by Cholesky decomposition and incorporate the correlation into the representation of the input random field $Y_j$ as follows \cite{he2016stochastic}
\begin{equation}
\label{KLE}
Y_j(x,\omega)=\bb{E}[Y_j(x,\omega)]+\sigma_{Y_j}\bigg(\sum_{k=1}^jL_{j,k}Y'_k(x,\omega)\bigg).
\end{equation}
It is noted that the KLE has the same mean and covariance function as the prescribed input field $Y_j$. Furthermore, the input random field $Y_j$ can be approximated by truncating the infinite summation in (\ref{KLE}) at the $N_j-$th term,
\[
Y_j(x,\omega)=\mathbb{E}[Y_j(x,\omega)]+\sigma_{Y_j}\bigg(\sum_{k=1}^jL_{j,k}\sum_{i=1}^{N_k} \sqrt{\lambda_i^k}e_i^k(x)Z^k_i(\omega)\bigg).
\]
As a result, $k(x,\omega)$ and $q(x,\omega)$ can be approximated by the representations with finite  parameters, i.e., $q(x,\omega)=\exp(Y_1(x,Z(\omega)))$ and $k(x,\omega)=\exp(Y_1(x,Z(\omega)))$, where $Z=(Z_1^1,\cdots,Z_{N_1}^1, Z_1^2,\cdots,Z_{N_2}^2)$. Thus, the stochastic property of the input random fields $Y_j(x,\omega)$  ($j=1,2$) can be effectively approximated with finite random variables $Z_i^j$  ($i=1,2,\cdots,N_j$;   $j=1,2$) by retaining the leading KLE terms with large eigenvalues.

\subsection{GMsFEM}
\label{GMsFEM}

\begin{figure}[htbp]
  \label{coarse-cell}
  \centering
  \includegraphics[width=3.3in, height=2in]{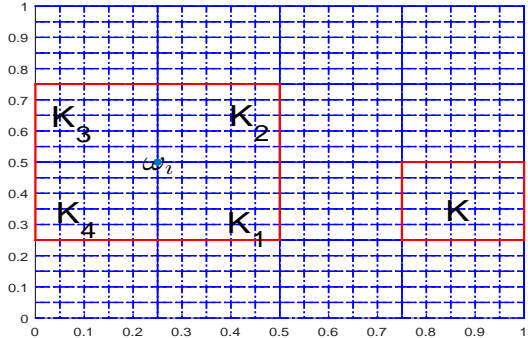}
  \caption{Illustration of a coarse neighborhood and a coarse element}
\end{figure}

 GMsFEM can achieve efficient forward model simulation and provide an accurate  approximation for  the solution of  multiscale problems.
  In this section, we follow the idea of  GMsFEM \cite{efendiev2013generalized} and apply it  to the time-fractional diffusion equation (\ref{frac-eq}).
  For GMsFEM, we need to pre-compute a set of multiscale basis functions. For the space $V^{\omega_i}_{\text{snap}}$ of snapshots,
  we solve the following local eigenvalue problem on each coarse block $\omega_i$,
\begin{eqnarray}
\label{MS-basis}
\begin{cases}
& -\text{div}(k(x, \mu_j)\nabla \varphi_{l,j}) =\lambda_{l,j} k(x, \mu_j) \varphi_{l,j}\ \ \text{in}\ \omega_i\\
& k(x, \mu_j)\nabla \varphi_{l,j}\cdot \vec{n} = 0\ \ \text{on}\ \partial \omega_i,
\end{cases}
\end{eqnarray}
where the samplers $\{\mu_j\}^{N_\mu}_j$ are drawn from the prior distribution of $\mu$. After the discretization using finite element method,  the local eigenvalue problem can be formulated as an
algebraic system,
\[
A(\mu_j)\varphi_{l,j}=\lambda_{l,j} S(\mu_j)\varphi_{l,j},
\]
where
\[
A(\mu_j)=[a(\mu_j)_{mn}]=\int_{\omega_i} k(x, \mu_j)\nabla v_n\nabla v_m, \quad S(\mu_j)=[s(\mu_j)_{mn}]=\int_{\omega_i} k(x, \mu_j) v_nv_m,
\]
where $v_n$ denotes the basis functions in fine grid. We take the first $M^i_{\text{snap}}$ eigenfunctions corresponding to the dominant eigenvalues for each coarse neighborhood $\omega_i$, (see Figure \ref{coarse-cell}),
 $i=1,2,\cdots,N_H$,  where $N_H$ is the number of coarse nodes. Hence we can construct the space of snapshots as
\[
V^{\omega_i}_{\text{snap}}=\text{span}\{\varphi_{l,j}, 1 \leq j \leq N_\mu, 1 \leq l \leq M^i_{\text{snap}}\}.
\]
and the snapshot functions can be stacked into a matrix as
\[
R_\text{snap}=[\varphi_1,\cdots,\varphi_\text{Msnap}],
\]
where $M_{\text{snap}}=N_{\mu}\times M^i_{\text{snap}}$ denotes the total number of snapshots used in the construction.

Next we perform a spectral decomposition of the space of snapshots and solve the local problems
  \begin{eqnarray}
   \begin{cases}
   & -\text{div}(k(x, \bar{\mu})\nabla \psi^i_k) =\lambda_k k(x, \bar{\mu}) \psi^i_k\ \ \text{in}\ \omega_i\\
   & k(x, \bar{\mu})\nabla \psi^i_k\cdot \vec{n} = 0\ \ \text{on}\ \partial \omega_i.
  \end{cases}
  \end{eqnarray}
where $\bar{\mu}=\frac{1}{N_{\mu}}\sum_{j=1}^{N_{\mu}}\mu_j$, $\psi_k^i \in V^{\omega_i}_{\text{snap}}$, to reduce the dimension of the snapshot space and construct the space $V^{\omega_i}$. We choose the smallest $M_i$ eigenvalues of
  \[
  A\psi^i_k=\lambda_k S\psi^i_k,
  \]
and take the corresponding eigenvectors in the space of snapshots by setting $\psi_k^i=\sum_j \psi_{k,j}^i\varphi_j$, for $k=1,\cdots,M_i$, to form the reduced snapshot space, where $\psi_{k,j}^i$ are the coordinates of the vector $\psi^i_k$ and
  \[
  A=[a_{mn}]=\int_{\omega_i} k(x, \bar{\mu})\nabla \varphi_n\nabla \varphi_m=R_{\text{snap}}^T\bar{A}R_{\text{snap}},
  \]
  \[
  S=[s_{mn}]=\int_{\omega_i} k(x, \bar{\mu}) \varphi_n\varphi_m=R_{\text{snap}}^T\bar{S}R_{\text{snap}}.
  \]
where $\bar{A}$ and $\bar{S}$ denote fine-scale matrices corresponding to the stiffness and mass matrices with the permeability $k(x, \bar{\mu})$. For each coarse element $K\in\omega_i$, let $\chi_i$ be the solution to the equation
  \[
    \left\{
  \begin{aligned}
  -\text{div}(k\nabla \chi_i)&=0\ \ K\in\omega_i\\
  \chi_i&=g_i\ \ \text{on}\ \partial K,
  \end{aligned}
     \right.
  \]
where $g_i$ is a  linear hat function.  The relationship between a coarse neighborhood and its coarse elements is illustrated in  Figure \ref{coarse-cell}.
 Thus $\{\chi_i\}_{i=1}^{N_H}$ form a set of partition of unity functions associated with the open cover $\{\omega_i\}_{i=1}^{N_H}$ of $\Omega$.
Then we  multiply the partition of unity functions by the eigenfunctions to construct GMsFE space,
  \[
  V_H=\text{span}\{\Psi_{il}: \Psi_{il}=\chi_i\psi^i_l: 1\leq i\leq N_H \quad \text{and}\quad 1\leq l\leq M_i\}.
  \]
We use  a single index  for  the multiscale basis function set  $\{\Psi_{il}\}$ and place them in  the following matrix
  \[
  R=[\Psi_1, \Psi_2,\cdots,\Psi_{M_v}],
  \]
where $M_v=\sum_{i=1}^{N_H} M_i$ denotes the total number of multiscale basis  functions.
 We note that once the matrix $R$ is constructed, it can be repeatedly used for simulation.

In the paper, we first discretize  the fractional  derivative  $(\ref{cp-div})$ as
 \begin{equation*}
  \begin{split}
   \int_0^t\frac{\partial u(x, s)}{\partial s}\frac{ds}{(t-s)^\gamma}&=\sum_{k=1}^n\int_{t_k}
   ^{t_{k+1}}\frac{\partial u(x, s)}{\partial s} \frac{ds}{(t-s)^\gamma} \\
   &=\sum_{k=1}^n \frac{u^{k+1}-u^k}{\Delta t} \int_{t_k}^{t_{k+1}} \frac{ds}{(t_{n+1}-s)^\gamma}+O(\Delta t) \\
   &=\frac{{\Delta t}^{-\gamma}}{1-\gamma}\sum_{k=1}^n(u^{k+1}-u^k)[(n+1-k)^{1-\gamma}-(n-k)^{1-\gamma}]+O(\Delta t),
  \end{split}
 \end{equation*}
where $0=t_1<t_2\cdots<t_{n+1}=t$ and $t_k=(k-1)\Delta t$. Thus,  we have the approximation of the derivative at time $t_{n+1}$ as
 \[
 ^c{D_t^\gamma} u (x, t_{n+1})=\frac{{\Delta t}^{-\gamma}}{\Gamma(2-\gamma)}\sum_{k=1}^n(u^{k+1}-u^k)[(n+1-k)^{1-\gamma}-(n-k)^{1-\gamma}]+O(\Delta t).
 \]
For simplicity of presentation, we define
 \[
 b_k:=(n+1-k)^{1-\gamma}-(n-k)^{1-\gamma}, \ \ k=1,2,\cdots,n.
 \]
Let $U^n$ be  the solution  at the $n-$th time level. Then we have the weak formulation for the
diffusion equation (\ref{frac-eq}),
  \[
    \left\{
  \begin{aligned}
  \frac{\Delta t^{-\gamma}}{\Gamma(2-\gamma)}\tilde{a}\bigg(\sum_{k=1}^n(U^{k+1}-U^k)b_k, v\bigg)+a(U^{n+1}, v)&=(f(t_n), v),\quad \forall v\in V_H\\
  (U^0,v)&=(u(x,0), v),\quad \forall v\in V_H,
  \end{aligned}
     \right.
  \]
where $(\cdot , \cdot)$ denotes the usual $L_2$ inner product and
  \[
    a(u,v)=\int k(x)\nabla u\nabla vdx, \ \
    \tilde{a}(u,v)=\int q(x) uvdx.
  \]
Let $c_k:=b_k-b_{k-1}$ and $b_0=0$. Then the weak formulation can be rewritten as
  \[
  \tilde{a}(U^{n+1}, v)+\Delta t^\gamma \Gamma(2-\gamma)a(U^{n+1}, v)=\tilde{a}\big(\sum_{k=1}^nU^kc_k, v\big)+\Delta t^\gamma \Gamma(2-\gamma)\big((f(t_{n+1}), v)\big).
  \]
We assume that $U^n$ has the expansion
  \[
  U^n=\sum_{j=1}^{M_v} u_{Hj}^n \Psi_j(x),
  \]
where the subscript $H$ denotes the GMsFEM solution on coarse grid. Let
  \[
  u_H^n=(u_{H1}^n,u_{H2}^n,\cdots,u_{HM_v}^n)^T.
  \]
Then for $k=1,\cdots,M_v$,
  \begin{eqnarray}
  \label{c-eq}
  \begin{split}
  &\sum_{j=1}^{M_v} u_{Hj}^{n+1}\tilde a(\Psi_j, \Psi_k)+\Delta t^\gamma \Gamma(2-\gamma)\sum_{j=1}^{M_v}u_{Hj}^n a(\Psi_j, \Psi_k)\\
  &=\sum_{i=1}^n\sum_{j=1}^{M_v} c_i u_{Hj}^{i} \tilde a(\Psi_j, \Psi_k)+\Delta t^\gamma \Gamma(2-\gamma)(f^{n+1}, \Psi_k).
  \end{split}
  \end{eqnarray}
 Let  $B$, $K$ and $F$ be the weighted mass, stiffness matrices and load vector using FEM basis functions in fine grid,  respectively. Then equation (\ref{c-eq}) gives
   the following algebraic system,
  \[
   R^TBRu_H^{n+1}+\Delta t^\gamma \Gamma(2-\gamma)R^TKRu_H^{n+1}=\sum_{i=1}^nc_iR^TBRu_H^{i}+\Delta t^\gamma \Gamma(2-\gamma)R^TF.
  \]
 If we denote
 \[\tilde B=R^TBR,  \quad \quad \tilde K=R^TKR,\]
  then $u_H^n$ can be computed  by the iteration
  \begin{equation}
  \label{iteration_c}
  u_H^{n+1}=\bigg(\tilde B+\Delta t^\gamma \Gamma(2-\gamma)\tilde K\bigg)^{-1}\bigg(\sum_{i=1}^nc_i\tilde{B}u_H^{i}+\Delta t^\gamma \Gamma(2-\gamma) R^TF\bigg).
  \end{equation}
By using the multiscale basis functions,  the solution in fine grid can be downscaled  by the transformation $Ru_H^n$.

We note that when GMsFEM is not applied, the full model solution is obtained by the iteration
  \begin{equation}
  \label{iteration_f}
  u_h^{n+1}=\bigg(B+\Delta t^\gamma \Gamma(2-\gamma)K\bigg)^{-1}\bigg(\sum_{i=1}^nc_iBu_h^{i}+\Delta t^\gamma \Gamma(2-\gamma)F\bigg).
  \end{equation}
Compared  $(\ref{iteration_c})$ with $(\ref{iteration_f})$, it can be seen that  the size of $\tilde K$ and $\tilde B$ are $M_v \times M_v$, but the size of $K$ and $B$ are $N_h \times N_h$  ($M_v\ll N_h$).
Thus a much smaller system is solved in GMsFEM.  The matrix $R$ for multiscale basis functions is computed overhead and it can be repeatedly used for all the time levels. This significantly improves the efficiency for forward model simulations.

\subsection{Stochastic collocation via least-squares method}
\label{SCM}
 Stochastic collocation method is an efficient approach  to approximate the solution of PDEs with random inputs. In the paper, the stochastic collocation method
 is based on generalized polynomial chaos (gPC) and least-squares method.   The approximation solution  can be represented by gPC expansion using the stochastic collocation method.
 We use the stochastic collocation method to solve the forward model.
 With the established gPC expansion of the approximation of the forward model,  the evaluation of the likelihood function $\bff{L}(z)$ in  MCMC sampling can be significantly accelerated.

We denote the random  parameters as $Z=(Z_1, \cdots, Z_{n_z})$, and assume that each random variable $Z_i$ has a prior probability density function $\pi_i(z_i): \Gamma_i \rightarrow \mathbb{R}$, for $i=1,\cdots, n_z$, where $\Gamma_i$ is the support of $Z_i$. Then the joint prior density function of $Z$ is
  \[
  \pi(z)=\prod_{i=1}^{n_z} \pi_i(z_i),
  \]
and its support has the form
  \[
  \Gamma :=\prod_{i=1}^{n_z} \Gamma_i \in \bb R^{n_z}.
  \]
For standard normal random parameters, the support of the prior is $\bb R^{n_z}$, and Hermite polynomials can be used.
If the support random parameters are bounded, e.g., $\Gamma:=\prod_{i=1}^{n_z} [a_i, b_i]^{n_z}$, we can use Legendre orthogonal polynomials as the basis functions to construct approximations of the forward model solution since random variables with bounded support in $[a, b]$ can always be mapped to $[-1, 1]$.

Without loss of generality, we describe the gPC approximation to the forward model for $n_d=1$. Let $i=(i_1,\cdots,i_{n_z})$ $\in$ $N_0^{n_z}$ be a multi-index with $|i|=i_1+\cdots+i_{n_z}$, and let $N\geq0$ be an integer. The $N$th-degree gPC expansion of $\bff{G}(Z)$ is defined as
  \begin{equation}
  \label{gpc_t}
    \bff{G}_N(Z)=\sum_{i=1}^P \bff{c}_i \Phi_i(Z), \quad \quad P=\frac{(N+n_z)!}{N!n_z!},
  \end{equation}
where $\Phi_i(Z)$ are the basis functions defined as
  \[
  \Phi_i(Z)=\phi_{i_1}(Z_1)\cdots\phi_{i_{n_z}}(Z_{n_z}), \quad 0\leq|i|\leq N,
  \]
and satisfies
  \begin{equation}
  \label{ortho-basis}
  \mathbb{E}[\Phi_i(Z)\Phi_j(Z)]=\int \Phi_i(z)\Phi_j(z)\pi(z)dz=\delta_{i,j}, \quad 0\leq |i|,|j|\leq N,
  \end{equation}
where $\delta_{i,j}=\prod_{k=1}^{n_z} \delta_{{i_k},{j_k}}$ because  $\phi_m(Z_k)$ is the $m-$th degree one-dimensional orthogonal polynomial having been normalised in the $Z_k$ direction.

In the stochastic collocation method, we first choose a set of collocation nodes $\{z^{(i)}\}_{i=1}^Q\in\Gamma$, where $Q\geq1$ is the number of nodes. Then  for each $i=1,\cdots,Q$, we solve a deterministic problem at the node $z^{(i)}$ to obtain
  \[
  \bff{G}(z^{(i)})=\bff g\circ u(x,t;z^{(i)}),
  \]
where $\bff g: \mathbb{R}^{n_u}\rightarrow \mathbb{R}$ is a state  function. After the pairings $\big(z^{(i)}, \bff{G}(z^{(i)})\big)$ ($ i=1,\cdots,Q$) being obtained, we are able to construct a good approximation of $\bff{G}_N(z)$, such that $\bff{G}_N(z^{(i)})=\bff{G}(z^{(i)})$ for all $i=1,\cdots,Q$. Thus, we need to solve $Q$ deterministic problems.
In the paper we use least-squares method to obtain the coefficient $\bff{c}$ in $(\ref{gpc_t})$.

Let $\{z^{(i)}\}_{i=1}^Q$ be the set of i.i.d. samples for $Z$ and $\{\bff{G}(z^{(i)})\}_{i=1}^Q$ the corresponding realizations of the stochastic function $\bff{G}(Z)$.
Let
\[
\bff c=(\bff{c}_1,\cdots,\bff{c}_P)^T\in \bb{R}^P,  \quad  \bff b=\bigg(\bff{G}(z^{(1)}),\cdots,\bff{G}(z^{(Q)})\bigg)^T\in \mathbb{R}^Q.
\]
If we  set the condition $\bff{G}_N(z^{(i)})=\bff{G}(z^{(i)})$,   $i=1,\cdots,Q$, then  the following equation holds,
  \begin{equation}
  \label{ls}
  \bff{Ac}=\bff b,
  \end{equation}
where $\bff A\in \bb{R}^{Q\times P}$ is the matrix with the entries
  \[
  \bff A_{ij}=\Phi_j(z^{(i)}), \quad i=1,\cdots,Q,\quad j=1,\cdots,P.
  \]
Solve the system in the least squares sense, the approximation coefficient is obtained as
  \begin{equation}
  \label{coef_ls}
  \tilde {\bff c}=(\bff{A}^T\bff{A})^{-1}\bff{A}^T\bff{b}.
  \end{equation}
Thus, we construct the $N$th-order gPC approximation by
  \begin{equation}
  \label{gpc_ls}
     \bff{\tilde G}_N(Z)=\sum_{m=1}^P  \bff {\tilde c_i} \Phi_i(Z).
  \end{equation}

When constructing the gPC approximation, the forward model will be solved $Q$ times to obtain the sampling vector $\bff b$. It can be obviously observed that $Q>P$, as we increase the order of the gPC expansion to pursue accuracy of the approximation, the simulation times we need to solve the deterministic forward model will increase significantly. We efficiently solve the problem using GMsFEM, i.e., $\bff{G}(z^{(i)})$ for $i=1,\cdots, Q$ are obtained by  GMsFEM.

\section{Numerical examples}

In this section, we use the intermediate distribution  strategy to build a reduced computational model  for the equation (\ref{frac-eq}),  and recover the model's inputs using Bayesian framework.
In Subsection \ref{rec-k}, we recover the permeability $k(x)$ when the fractional derivative   $\gamma$ known or unknown.
In Subsection \ref{rec-ks}, we recover the permeability $k(x)$ and specific fields $q(x)$ simultaneously when they are correlated each other.
 In Subsection \ref{sec4-3}, we will  identify a permeability  field with high permeability layer and low permeability layer mixed in the spatial space.

For the numerical examples,  we consider a dimensionless  square domain $\Omega:=[0,1]\times[0,1]$ for space, and $(0, T]$ for time, and we set the initial function as $u(x,0)=0$. Measurement data are generated by using finite element method in a  fine grid with time step $\Delta t = 0.001$, and the measurement noise is set to be $\sigma=0.01$. For any given values of $z$, we solve the time fractional diffusion equation  using GMsFEM with time step $\Delta t=0.002$. The regularization parameter $\alpha$ is set be  $0.01$ in the regularizing Levenberg-Marquart algorithm. For all numerical examples, the  parameters in DREAM$_\text{ZS}$ algorithm are set as follows:
\[
\beta = 2,\quad h_1 = 0.01, \quad h_2=10^{-8}.
 \]
 The number of samplers in constructing the polynomial surrogate is set as $Q=2P$, where $P$ is the number of the gPC basis functions for approximation. We will compare the results obtained by using LS-SCM based on the original Gaussian prior with the one
  using  intermediate distribution. For the convenience of presentation, we refer the surrogate  based on the original Gauss prior to  prior-based surrogate,
   and refer the surrogate based on the intermediate distribution to intermediate-based surrogate.

\subsection{Identification of the permeability and fractional derivative}
\label{rec-k}

\begin{figure}
\centering
  \includegraphics[width=3in, height=2in]{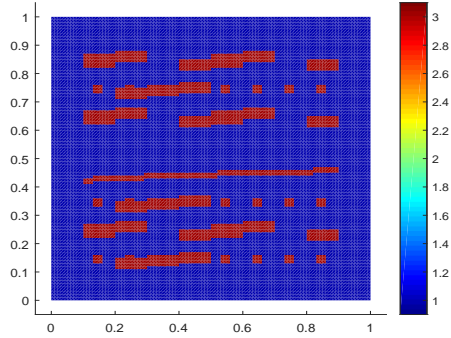}\\
  \caption{The spatial distribution of $\mathbb{E}[Y_1(x,\omega)]$}\label{mu}
\end{figure}

Consider the problem with mixed boundary condition, where  homogeneous Dirichlet boundary condition is imposed on $y=1$ or $x=1$, and the zero Neuman boundary condition is imposed on $y=0$ or $x=0$.  The specific field  $q(x)=1$ is given, source term is set as $f=10$,  the end time is $T=0.15$. For the first example, the fractional derivative is given by   $\gamma=0.5$.  The random field is generated with $l_{x,Y_1}=0.2$, $l_{y,Y_1}=0.4$, $\sigma_{Y_1}^2=1$, the mean of the random field is heterogeneous and shown as Figure \ref{mu}. The covariance function for the random field is given by
\begin{equation}\label{covfun}
\text{cov}[Y_1](x_1,y_1;x_2,y_2)=\sigma_{Y_1}^2\exp\bigg(-\frac{|x_1-x_2|^2}{2l_{x,Y_1}^2}-\frac{|y_1-y_2|^2}{2l_{y,Y_1}^2}\bigg),
\end{equation}
which is the Gaussian  kernel covariance.  We truncate the KLE by  the first $10$ terms,   i.e.,
\[
Y_1(x,\omega)=\log{k(x,\omega)}=\sum_{i=1}^{N_1} \sqrt{\lambda_i}e_i(x)z_i(\omega)+\mathbb{E}[Y_1(x,\omega)], \quad N_1=10.
\]
The ground true parameter $z_0$ are randomly drawn from the standard multivariate normal distribution. Measurements are taken at time instances  $0.02:0.01:0.11$ and the locations are distributed on the uniform $5\times 5$ grid of the domain $[0, 0.8]\times [0, 0.8]$.

The forward model is solved on a uniform $100\times100$ fine grid, and we set the coarse grid  $5\times5$ for GMsFEM simulation.  We set the number of training samplers $N_{\mu}=10$ and the number of eigenfunctions selected in calculating the snapshot space as $M^i_{\text{snap}}=20$. We select $M_c=5$ multiscale basis  functions at each block to construct the coarse reduced order model and solve the optimization problem (\ref{coarseopt}). The stepsize of input parameter is set as $h=0.5$ when we calculate the matrix $\bff{J}$ in the optimization problem
 by difference method.

We construct the intermediate distribution after obtaining the approximate sampling distribution, and we find that $2\sigma_\text{OLS}\sqrt\delta_k\leq1$ for $k=1,\cdots ,N_1$, the intermediate distribution is then set as $Z_k\sim U(z_\text{OLS$_k$}-1, z_\text{OLS$_k$}+1)$. Ten  local multiscale basis functions ( $M_i=10$)  are selected  in constructing the matrix $R$
to construct the gPC surrogate model. When the order of gPC is set as $N=3$, the number of random samplers is 572 when computing  vectors $\bff b$ and $\bff A$ in Section \ref{SCM}. We run 5 Markov chains simultaneously using DREAM$_\text{ZS}$ algorithm, and the length of each chain is $10000$, only the last 5000 realizations of each chain are used to compute the relevant statistical quantities.  We choose the  parameter $n_\text{CR}=8$ in  DREAM$_\text{ZS}$ for simulation.

\begin{figure}[htbp]
\centering
 \includegraphics[width=1.65in, height=1.4in]{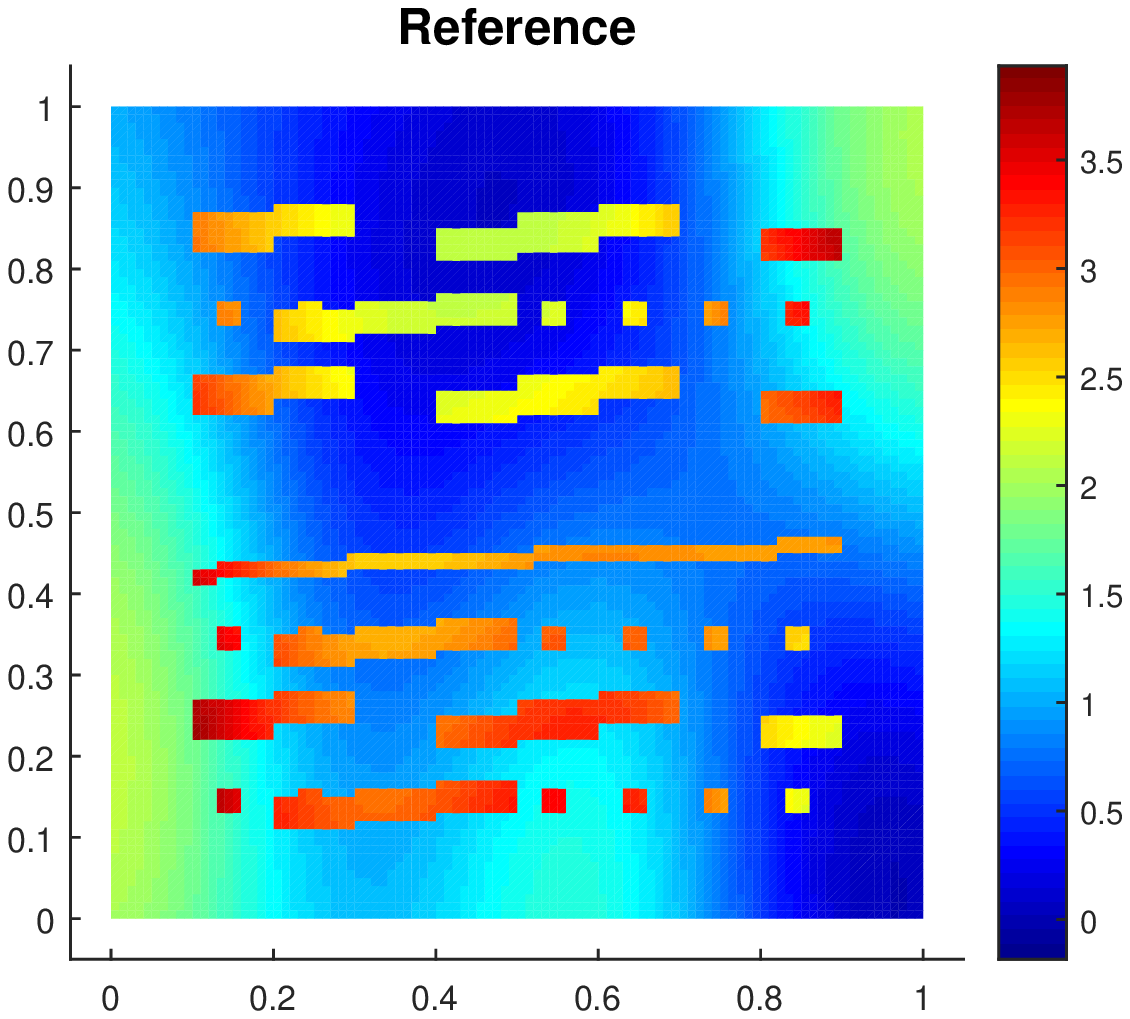}
  \includegraphics[width=1.65in, height=1.4in]{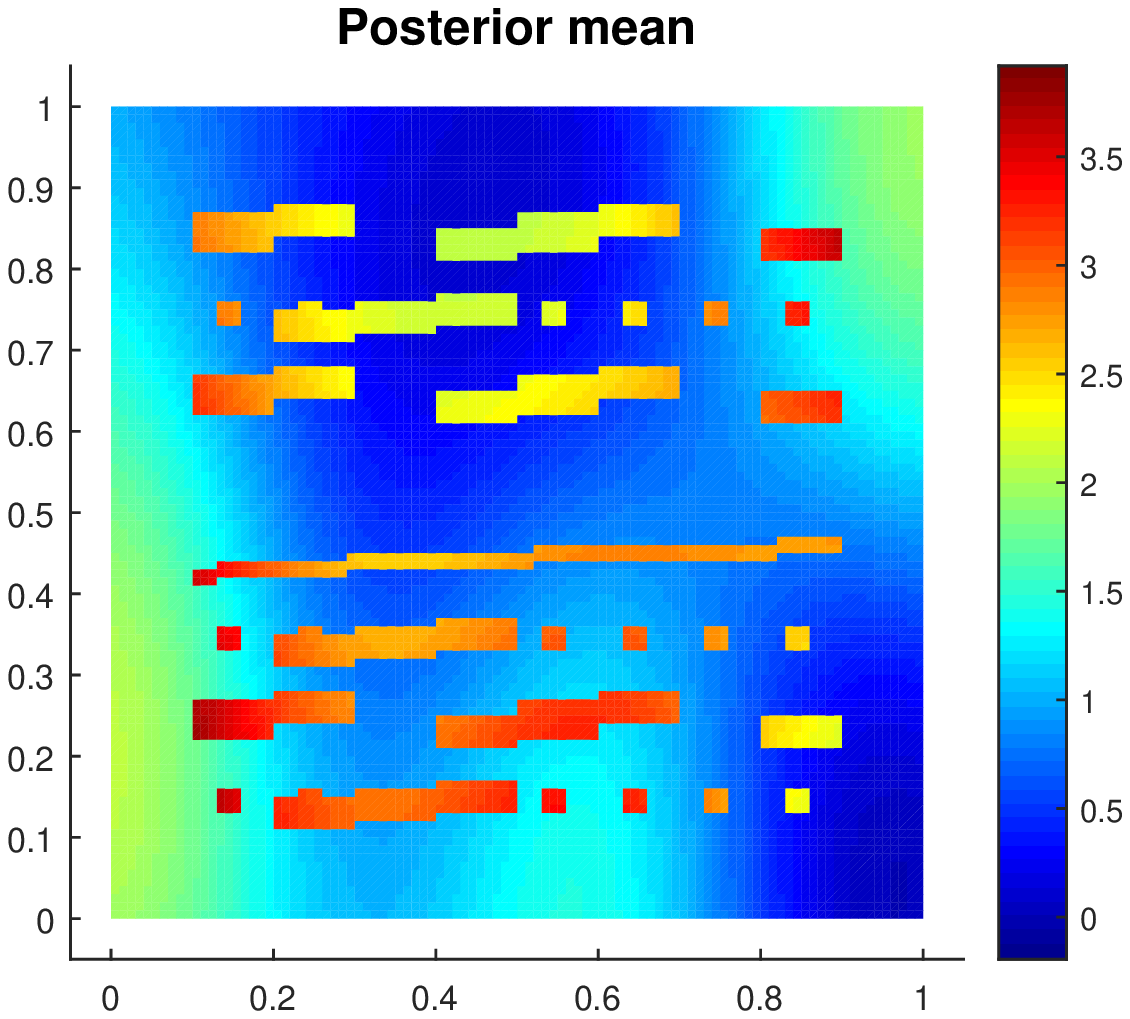}
  \includegraphics[width=1.65in, height=1.4in]{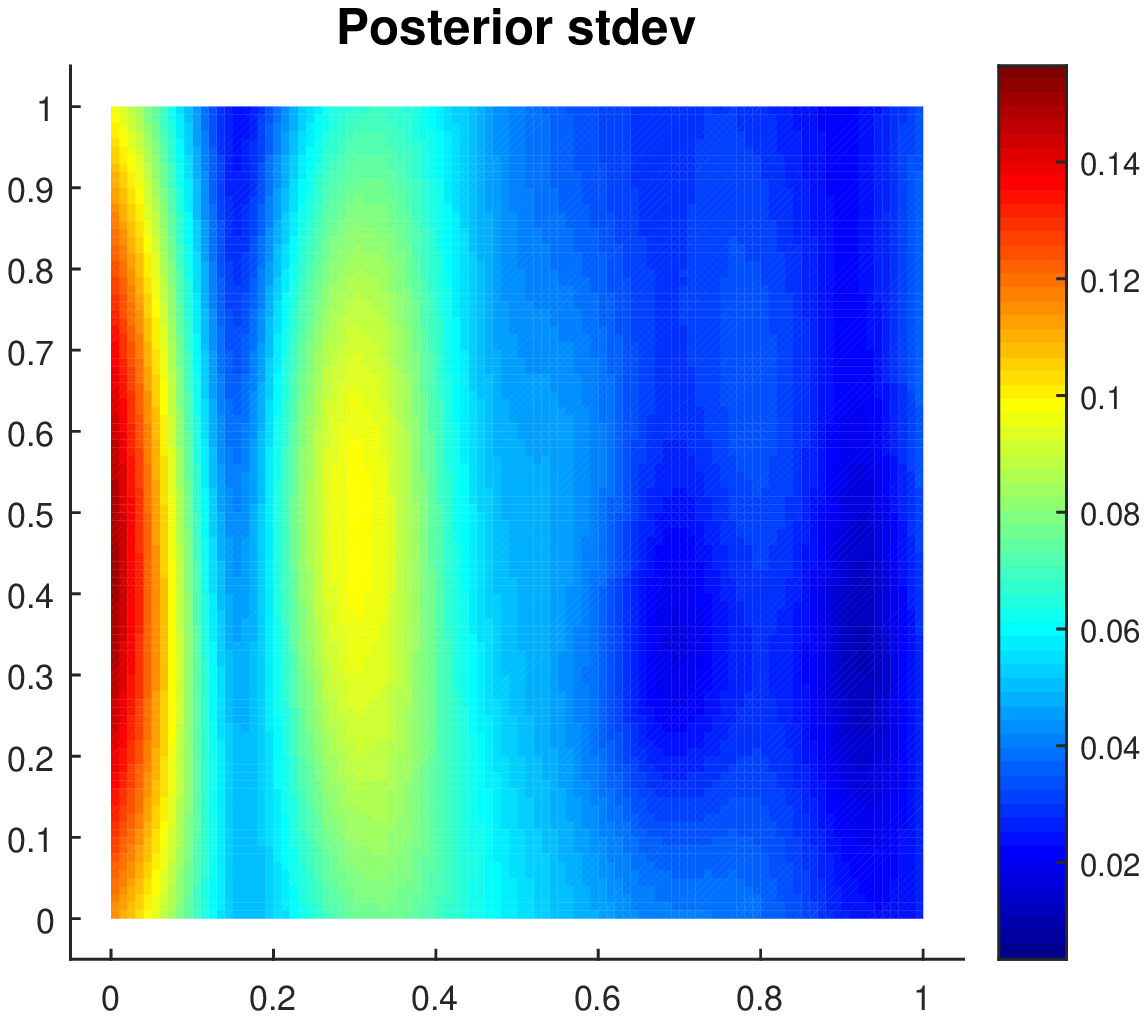}
\caption{True profile, posterior mean and posterior standard deviation of $\log k(x, \omega)$.}
  \label{coef1}
\end{figure}

\begin{figure}
\centering
  \includegraphics[width=0.97\textwidth, height=0.8\textwidth]{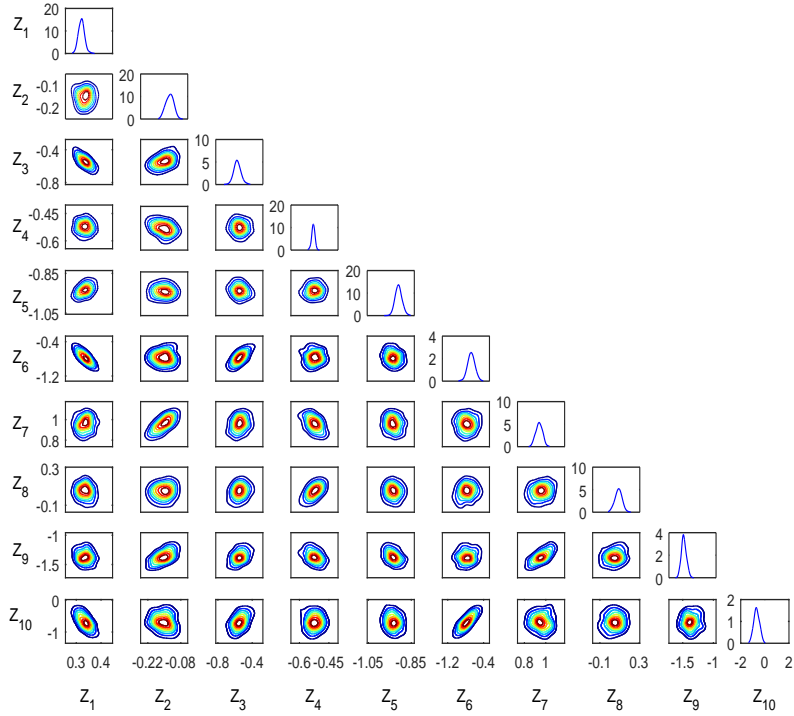}\\
  \caption{1-D and 2-D posterior marginals of $z$.}
  \label{biva}
\end{figure}

Figure \ref{coef1} shows the result using the proposed  strategy.  We compare MCMC mean estimate with the reference one in logarithmic scale and present the posterior standard deviation
in the figure.  We can see that the uncertainty is largest near the boundary $x=0$. Figure \ref{biva} shows all of the one and two- dimensional posterior marginals of $z$,  where some correlations are apparent between the lower-indexed modes and the higher-indexed modes, such as $z_3$ and $z_6$, $z_1$ and $z_{10}$, $z_2$ and $z_7$, but most of the modes appear uncorrelated and  mutually independent  based on the shape of their 2-D marginals.

For a comparison, we also construct the surrogate with respect to the original Gaussian prior distribution, i.e., the prior-based surrogate. Here we use a surrogate using order $N=3$ Hermite polynomials.   The same 572 training samplers are used for LS-SCM, and the model evaluations are implemented  by GMsFEM with $M_i=10$ local GMsFEM basis functions. Here we use
 the Gaussian prior
\[
\pi(z)\propto \exp(-\frac{z^Tz}{2}),
\]
and then we compute  the Kullback-Leibler (KL) divergence between the approximate posterior $\tilde \pi^d(z)$ and the exact one $\pi^d(z)$ for the surrogates.

\begin{figure}
\centering
  \includegraphics[width=0.6\textwidth]{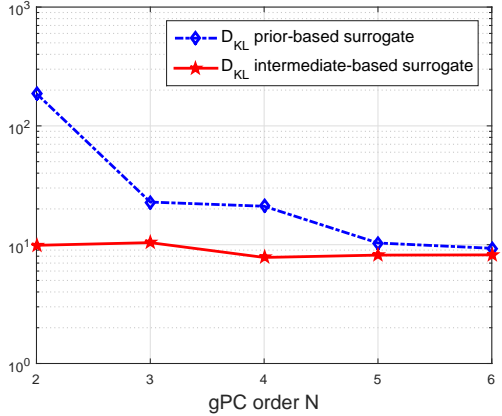}\\
  \caption{Kullback-Leibler divergence between the approximate posterior and the exact one, the blue dashed one is obtained by the prior-based surrogate, and the red solid one is obtained by the intermediate-based surrogate.}
  \label{KLD}
\end{figure}

For probability density functions $\pi_1(z)$ and $\pi_2(z)$, KL divergence is defined by
\[
D_{KL}(\pi_1||\pi_2)=\int \pi_1(z) \log \frac{\pi_1(z)}{\pi_2(z)}dz.
\]
We compute  the KL difference between $\tilde \pi^d(z)$ and $\pi^d(z)$ by
\[
D_{KL}(\tilde \pi^d||\pi^d)=\bb{E}_{\tilde\pi^d}[\log \frac{\bff{\tilde L}}{\bff L}]+\log \bb{E}_{ \tilde \pi^d}[\frac{\bff L}{\bff{\tilde L}}],
\]
where $\bff{\tilde L}$ refers to  the surrogate likelihood and $\bff L$ is the likelihood obtained via the full order model. The Kullback-Leibler divergence is plotted against the order of gPC in Figure \ref{KLD}. As a general trend from the figure,  the two curves of Kullback-Leibler divergence decrease as we increase the order of gPC. The main difference between them  is that Kullback-Leibler divergence  by the intermediate-based surrogate decreases slower than the one obtained by the prior-based surrogate.
 The KL divergence curves in Figure \ref{KLD} show that the posterior obtained by the intermediate-based surrogate approximates the reference  better than the one obtained by prior-based surrogate with the same gPC order.

The different performance of intermediate-based surrogate and prior-based surrogate may be caused by the data information  in the intermediate distribution. When we construct surrogate with respect to the intermediate distribution, we concentrate on the important region of the posterior, most part of the unimportant  region of posterior contained in the support of prior has been excluded when the intermediate distribution is constructed. For the example,  the average acceptance rate of MCMC chains for the intermediate-based surrogate is 15.04\%, and  6.70\% for the prior-based surrogate.  This also shows  the  advantage of the intermediate-based surrogate over  the prior-based surrogate.

It would be interesting to recover the fractional order $\gamma$ in the model (\ref{frac-eq}) simultaneously with the permeability field.  To this end, we keep the same conditions, but the unknown inputs are  $Z=[\gamma;Z_1,...,Z_{N_1}]$. The ground true value of $\gamma$ is 0.5, measurements are taken at the same locations and the same time levels as before. Since the unknown $\gamma\in (0, 1)$, we set the step size for the first parameter as $h=0.1$ and the step size for other components of $Z$ is still  $h=0.5$ for solving the optimization problem (\ref{coarseopt}).
 The coarse   model is constructed using $5$ local GMsFE basis functions (i.e., $M_c=5$) and the regularizing Levenberg-Marquart algorithm is used to  to solve the optimization problem. The motivation of our method is to construct surrogate on the important region of the posterior compared to the wide prior, we use $U(0,1)$ as the intermediate density for the first component because $(\gamma_\text{OLS}-1, \gamma_\text{OLS}+1)$ may be wider than the original  interval of $\gamma$. Hence, the intermediate distribution is set as $U(0,1)\otimes \prod_{i=1}^{N_1}U(z_\text{OLS$_i$}-1, z_\text{OLS$_i$}+1)$, and then we use Legendre orthogonal polynomials as the basis functions to construct the LS-SCM approximations of the forward model.

\begin{figure}
\centering
  \includegraphics[width=\textwidth, height=0.8\textwidth]{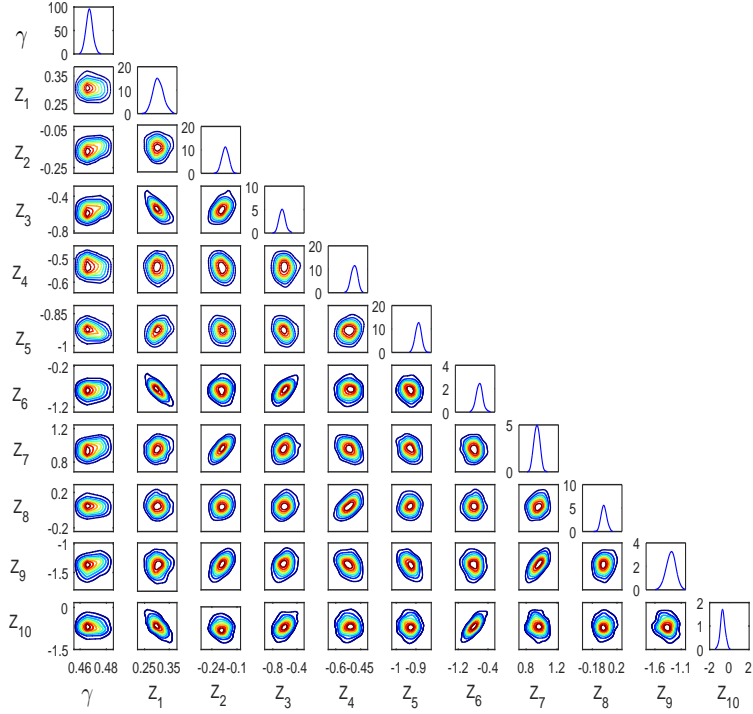}\\
  \caption{1-D and 2-D posterior marginals of $[\gamma;z]$.}
  \label{xbiva}
\end{figure}

In the following, we set the gPC order $N=3$ and $10$ ($M_i=10$) number of local multiscale basis functions  to construct the surrogate with respect to the intermediate distribution. We run 5 Markov chains using the surrogate model and remain the last 50\% of each chain for statistics computation.   Figure \ref{xbiva} illustrates the 1D and 2D marginal density of the parameters. It can be seen that except for very small correlation with $z_2, z_3, z_9$, the $\gamma$ is almost independent of all the parameters  in the permeability field.   The correlations between $zs$ are the same as the case when $\gamma$ is known compared to Figure \ref{biva}.

If the intermediate-based surrogate is not implemented, the surrogate can be constructed by using the Hermite orthogonal polynomials, i.e., when constructing the marginal matrix $\bff{A}$, samplers are taken from the standard Gaussian distribution, the range of component of the random vector  belongs to $\mathbb{R}$. Since $\gamma\in (0, 1)$, we make a transform \cite{marzouk2009dimensionality} for the first component of samplers, and we can choose the transform to be
\[
\bff{h}(x)=\frac{1}{2}+\frac{1}{\pi}\arctan(x),\ \ x\in\mathbb{R}.
\]
The $\bff{h}(x)$ is a bijection map and its range is $(0, 1)$.

\begin{figure}[htbp]
\centering
 \includegraphics[width=2in, height=1.7in]{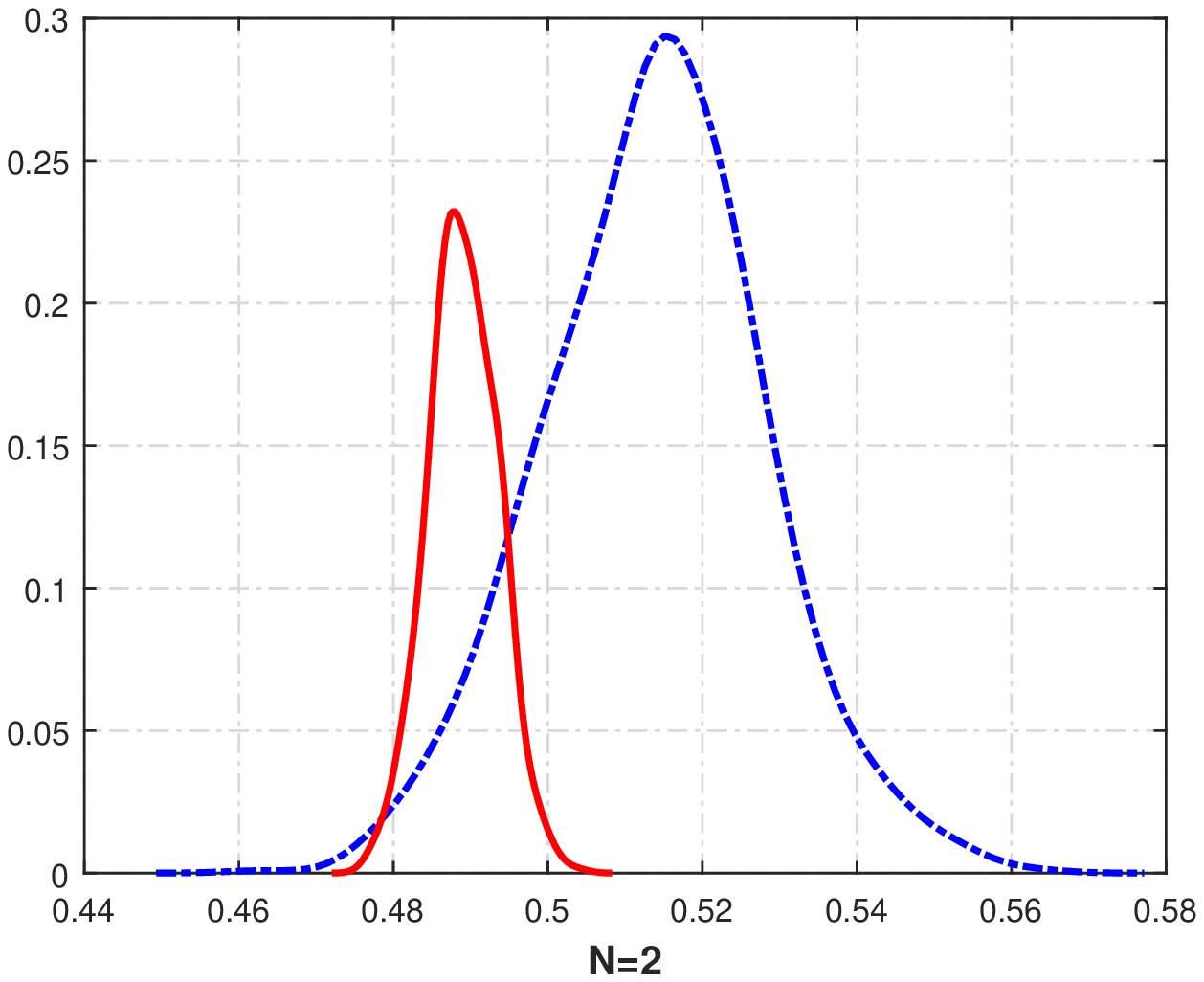}
  \includegraphics[width=2in, height=1.7in]{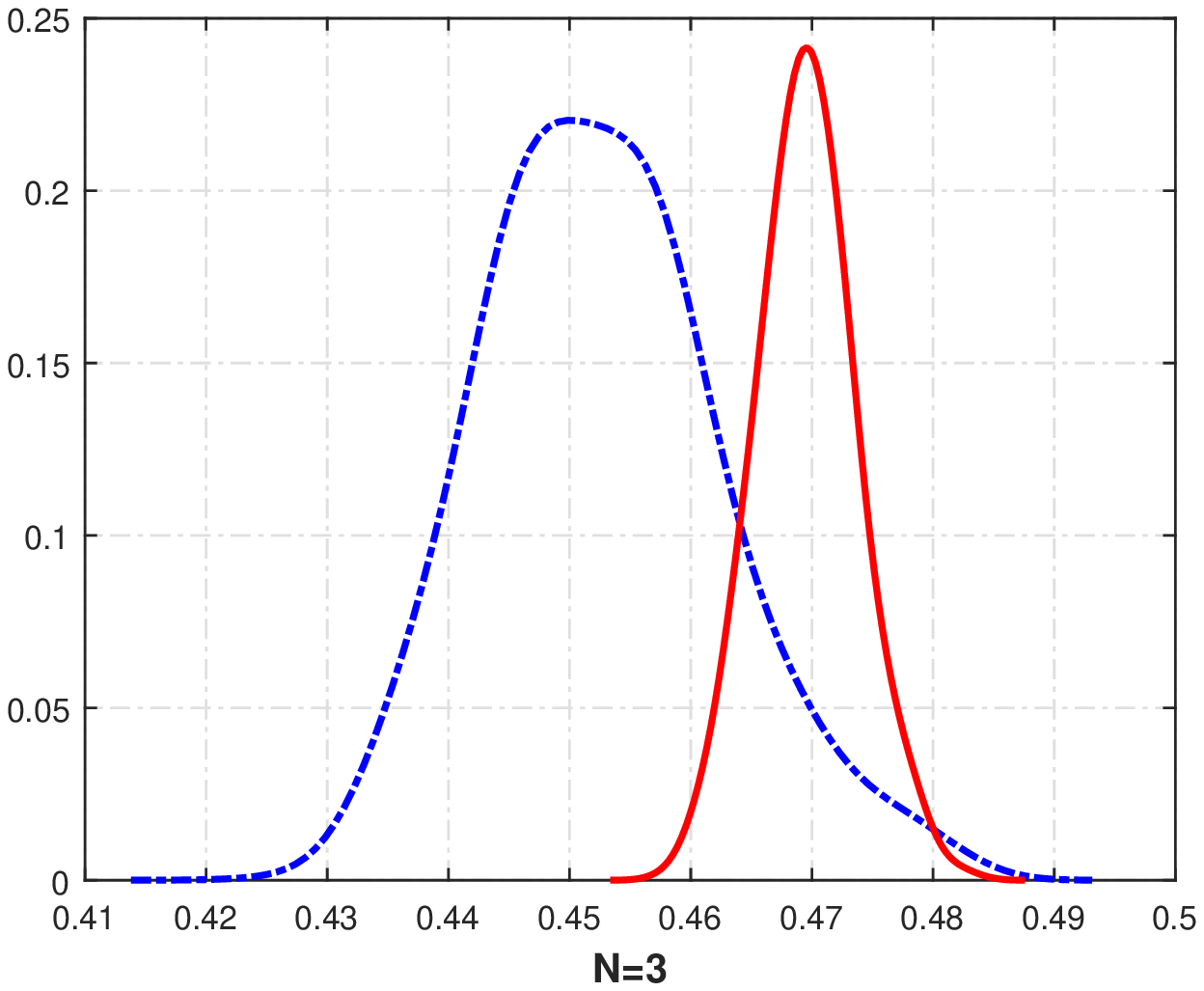}\\
  \includegraphics[width=2in, height=1.7in]{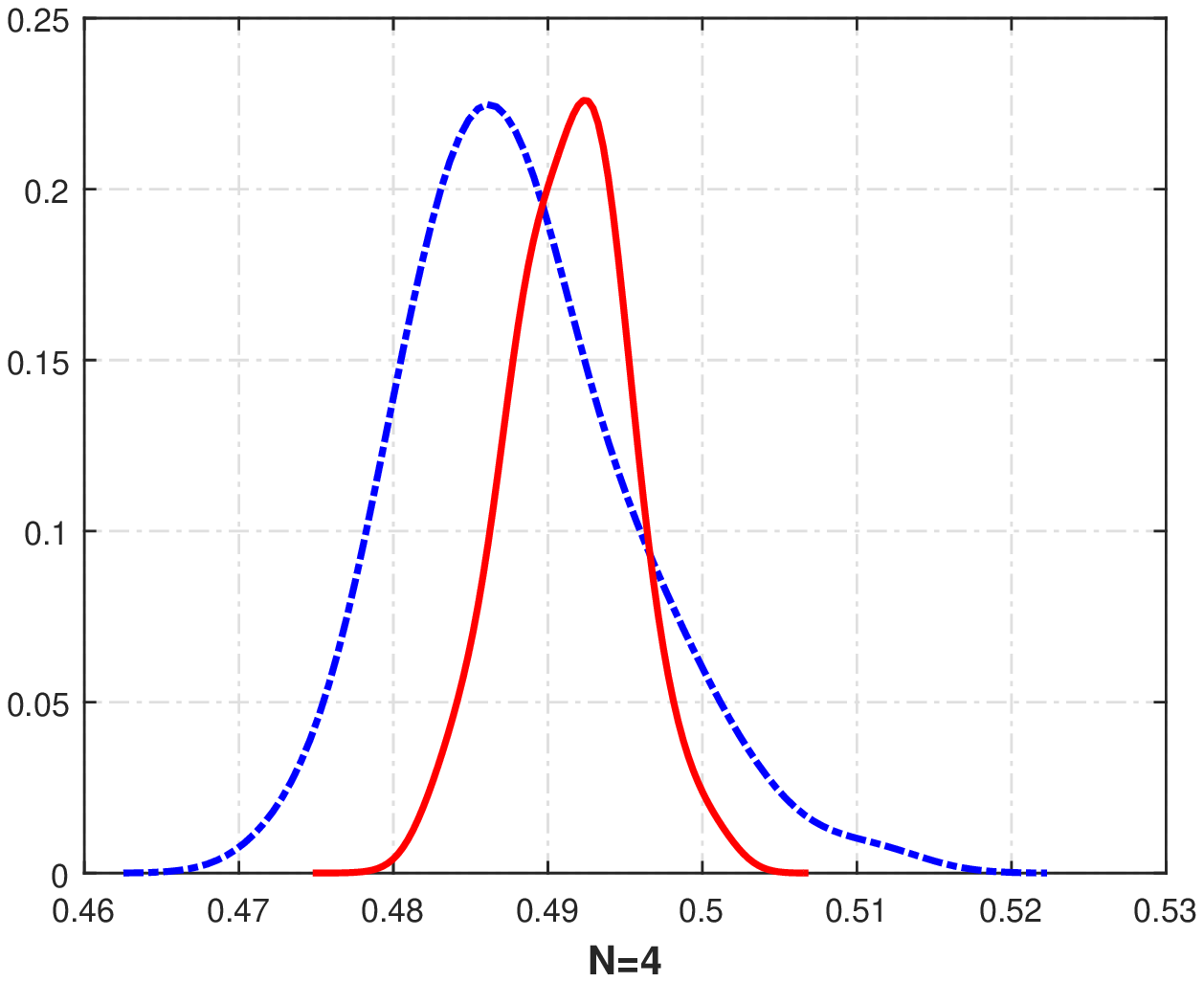}
  \includegraphics[width=2in, height=1.7in]{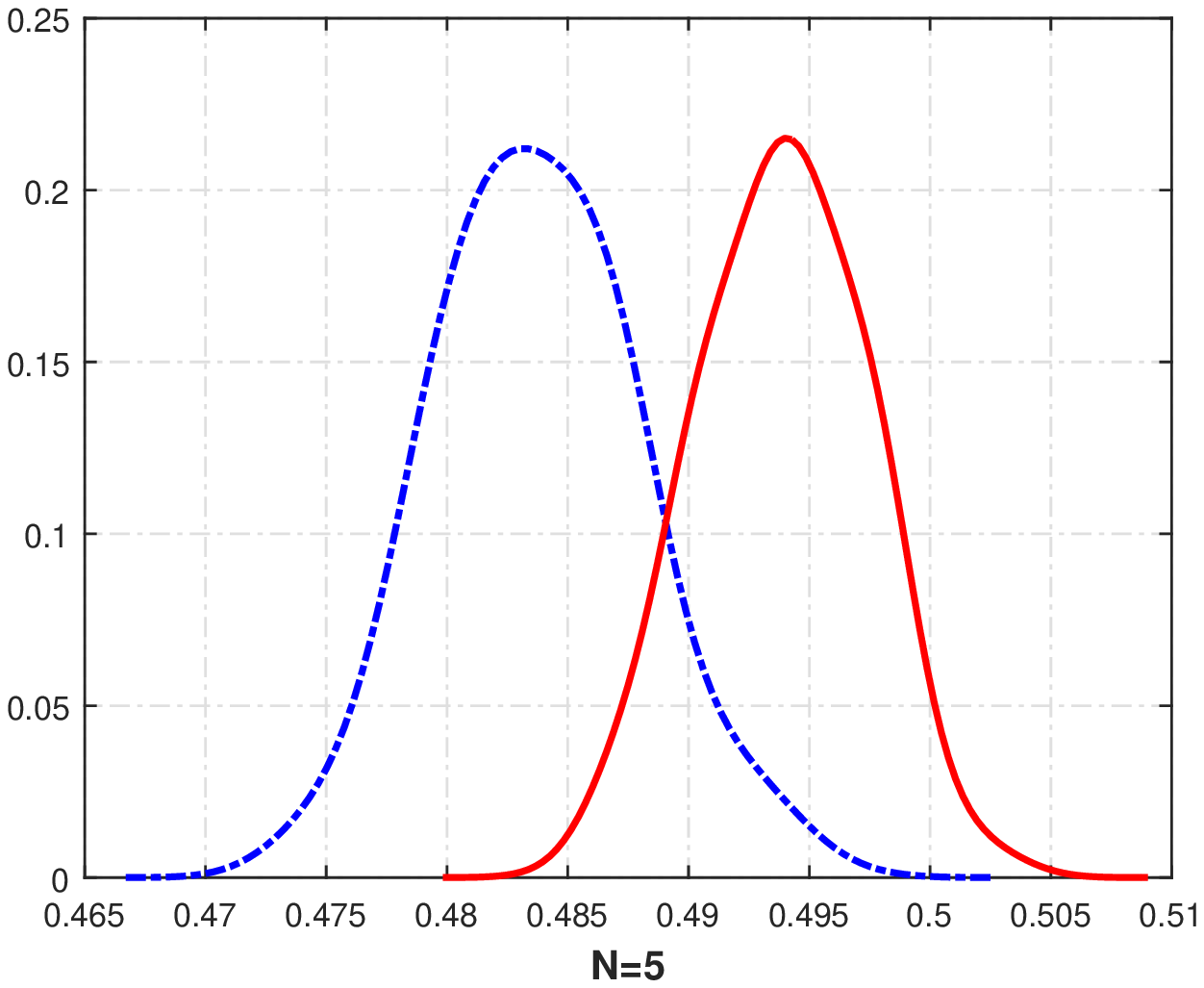}
  \caption{Marginal posterior density estimation of $\gamma$ with different gPC order, solid line is obtained via the intermediate-based surrogate and the dashed line is obtained via the prior-based surrogate.}
  \label{gamma}
\end{figure}

As shown in Figure \ref{gamma}, we consider the difference for  $\gamma$ first. We plot the posterior marginal density of $\gamma$ against the gPC order for the two types of surrogate, where the dashed one is computed  by prior-based surrogate, and the solid one is obtained through intermediate-based surrogate, and we denote them as $\tilde \pi^d_\text{pri}(\gamma)$ and $\tilde \pi^d_\text{int}(\gamma)$, respectively. It can be seen the support of $\tilde\pi^d_\text{int}(\gamma)$ is narrower than $\tilde \pi^d_\text{pri}(\gamma)$'s with different gPC order $N$. This implies $\gamma$ derived by the prior-based surrogate is more uncertain than the one derived by the intermediate-based surrogate, and we have more confidence in $\gamma$  by the intermediate-based surrogate.

\begin{figure}[htbp]
  \centering
  \includegraphics[width=0.6\textwidth]{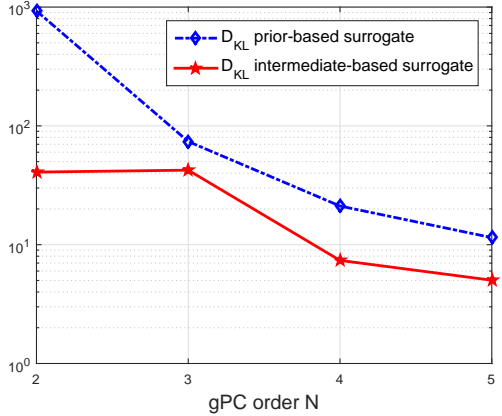}\\
  \caption{Kullback-Leibler divergence between the approximate posterior and the exact one, the solid one is based on the intermediate-based surrogate and the dashed one is computed  based on the prior-based surrogate}\label{KLD2}
\end{figure}

We also compute the KL divergence  between the reference posterior density and the approximate posterior density. The KL divergence is plotted against the gPC order in Figure \ref{KLD2}. It is obvious that both the two  approximate posteriors get  better  as the increase of gPC order.   With the same gPC order, the approximate posterior
 by the intermediate-based surrogate performs better than the prior-based surrogate. The proposed  strategy also has
 higher average acceptance rate of the 5 chains than the prior-based surrogate, this is due to the fact that we have rejected samplers having little occurrence probability through the approximate optimization problem, and we only focus on the small region where samplers have  high posterior probability.

\subsection{Identification of permeability field and specific field}
\label{rec-ks}
\begin{figure}
\centering
 \includegraphics[width=2in, height=1.7in]{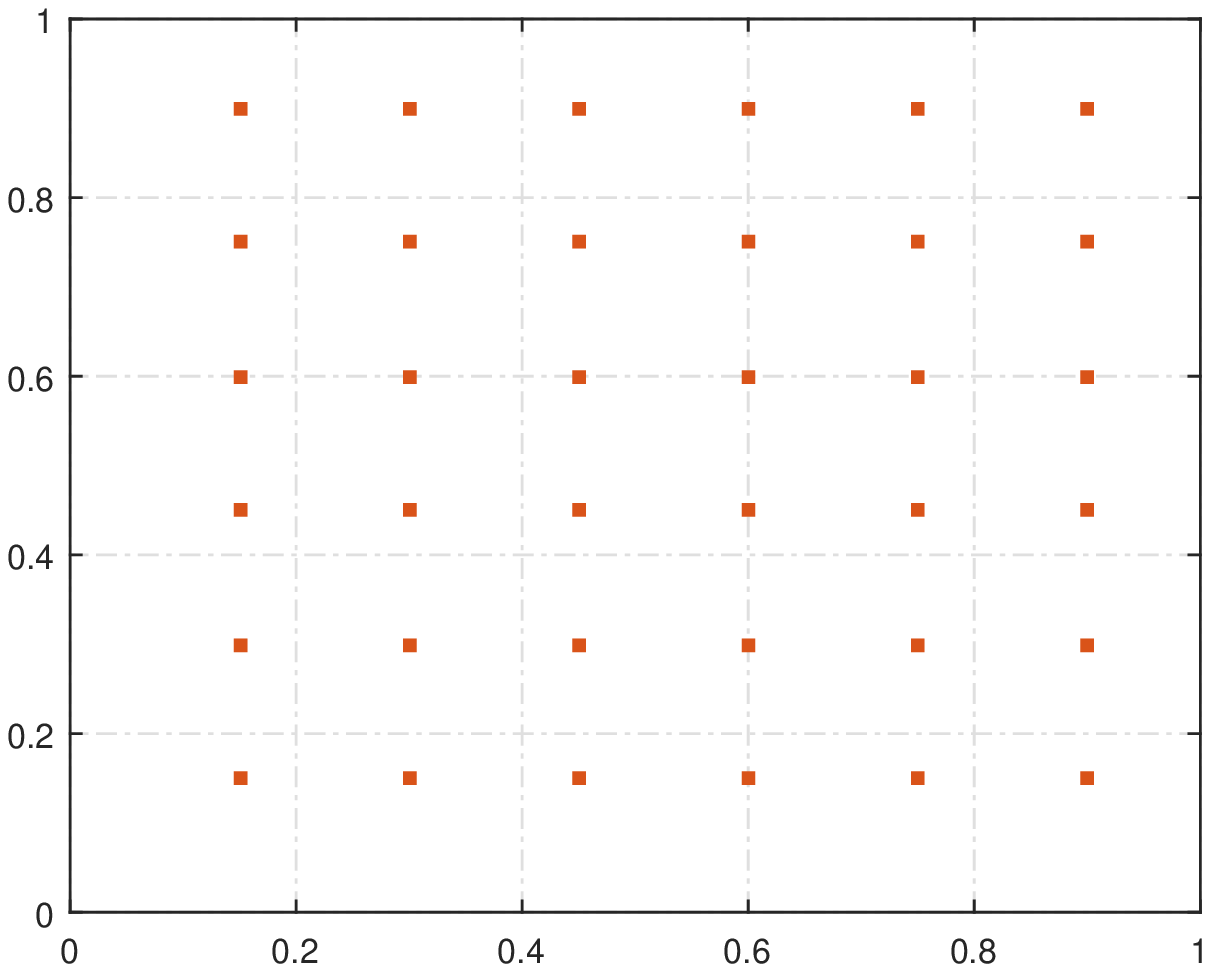}
  \includegraphics[width=2in, height=1.7in]{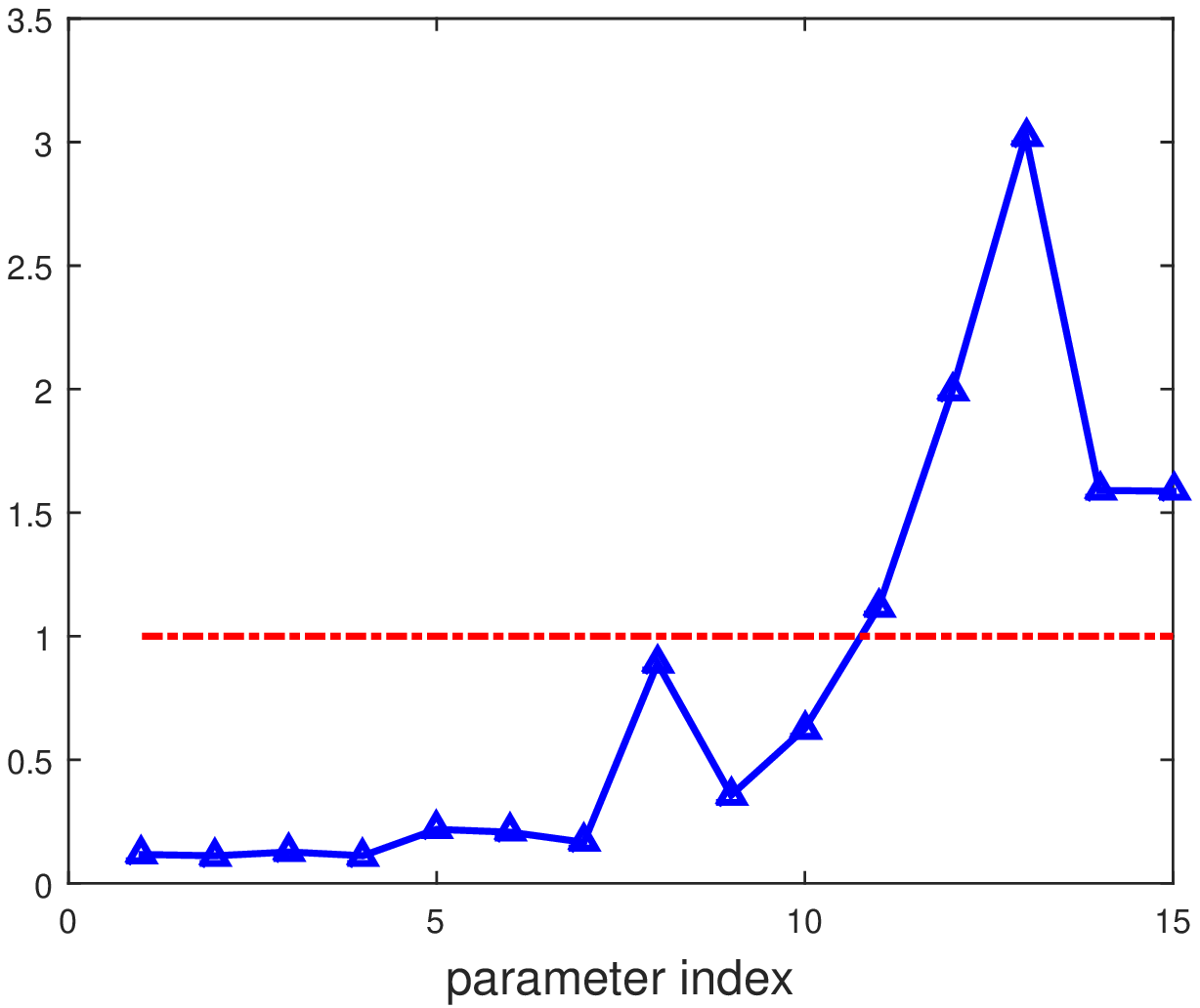}\\
\caption{Measurement locations for Subsection \ref{rec-ks} (left) and $2\sigma_\text{OLS}\sqrt\delta$ (right), the horizontal line is the criterion 1}\label{measure}
\end{figure}

In this subsection,  we consider the problem with homogeneous Dirichlet boundary condition and both the permeability and the specific fields are unknown. For the permeability field , we set $l_{x,Y_1}=0.1$, $l_{y,Y_1}=0.01$, $\sigma_{Y_1}^2=2$, $\bb{E}[Y_1(x,\omega)]=0$, and the covariance function is assumed to be
  \[
  \text{cov}[Y_1](x_1,y_1;x_2,y_2)=\sigma_{Y_1}^2\exp\bigg(-\frac{|x_1-x_2|}{2l_{x,Y_1}^2}-\frac{|y_1-y_2|}{2l_{y,Y_1}^2}\bigg),
  \]
which is the exponential kernel covariance.  For the specific field, we assume it to be smooth  and the covariance function has the Gaussian kernel   defined in Subsection \ref{rec-k}. The preset parameters are  $l_{x,Y_2}=0.1$, $l_{y,Y_2}=0.3$, $\sigma_{Y_2}^2=1$, $\mathbb{E}[Y_2(x,\omega)]=1$. We assume the two random fields are correlated with each other, and the correlation coefficient $\rho$ between them is -0.4. Their  truncated KLEs have the following form,
  \begin{eqnarray*}
  Y_1(x,\omega) &=& \log{k(x,\omega)}=\sigma_{Y_1}\bigg(\rho\sum_{i=1}^{N_2} \sqrt{\lambda^2_i}e^2_i(x)z^2_i(\omega)+
  \sqrt{1-\rho^2}\sum_{i=1}^{N_1} \sqrt{\lambda^1_i}e^1_i(x)z^1_i(\omega)\bigg), \\
  Y_2(x,\omega) &=& \log{q(x,\omega)}=\mathbb{E}[Y_2(x,\omega)]+\sigma_{Y_2}\sum_{i=1}^{N_2} \sqrt{\lambda^2_i}e^2_i(x)z^2_i(\omega).
  \end{eqnarray*}
Thus we have the unknowns $z=(z_1^1,...,z_{N_1}^1, z_1^2,...,z_{N_2}^2)$ with dimension $N_1+N_2$, we can see there are $N_1+N_2$ unknowns correlated with the permeability field and $N_2$ unknowns in the specific field from the expression of the two truncated random fields. Here we take  $N_1=7,N_2=8$ for simulation.   Measurements are taken at time instances $0.02 : 0.01 : 0.15$, the measurement locations are shown as Figure \ref{measure} (left),  which are measured from the boundary $x=0$ and $y=0$, with step size $30/200$.

\begin{figure}
  \centering
  \includegraphics[width=1.6in, height=1.7in]{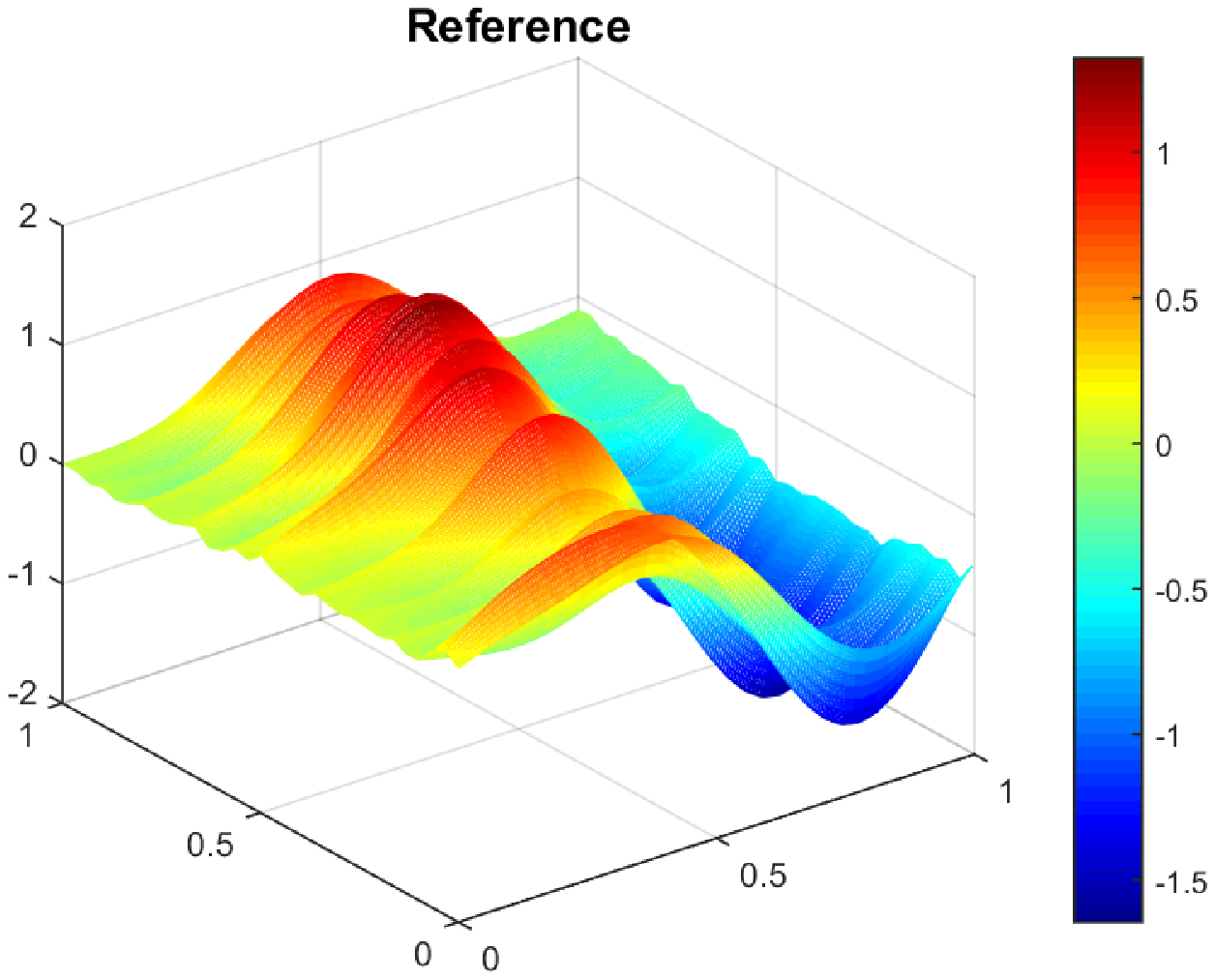}
  \includegraphics[width=1.6in, height=1.7in]{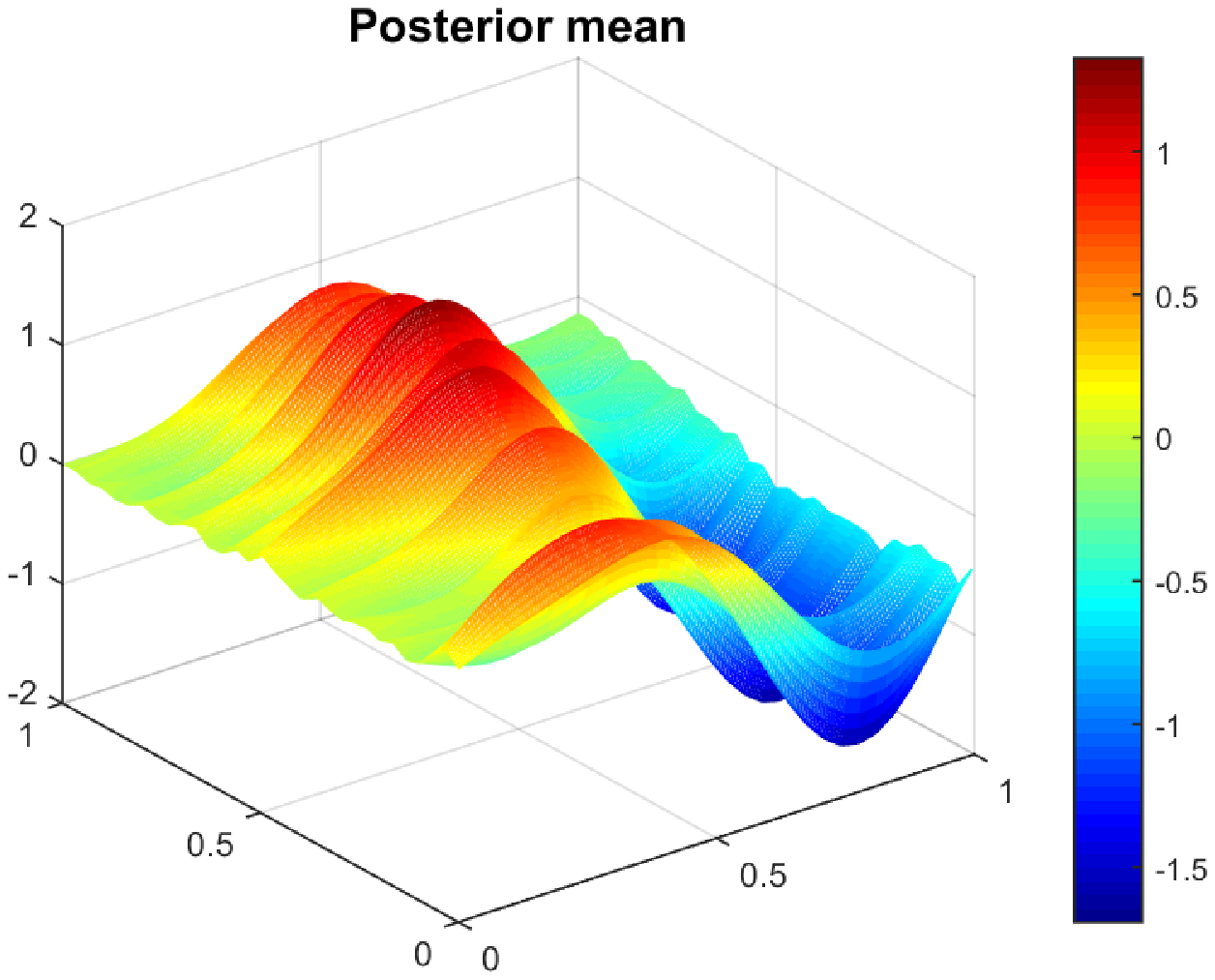}
  \includegraphics[width=1.6in, height=1.7in]{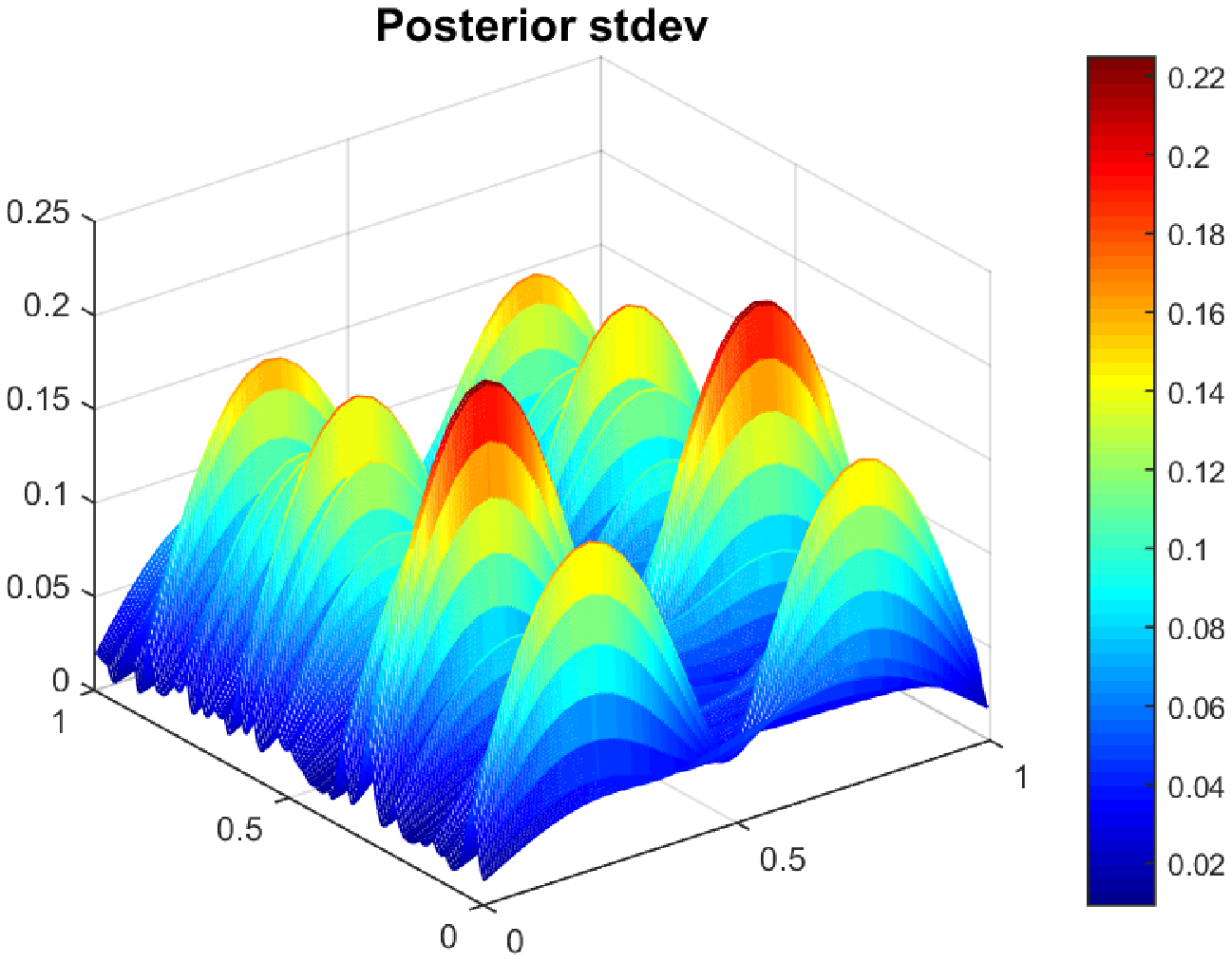}\\
  \includegraphics[width=1.6in, height=1.7in]{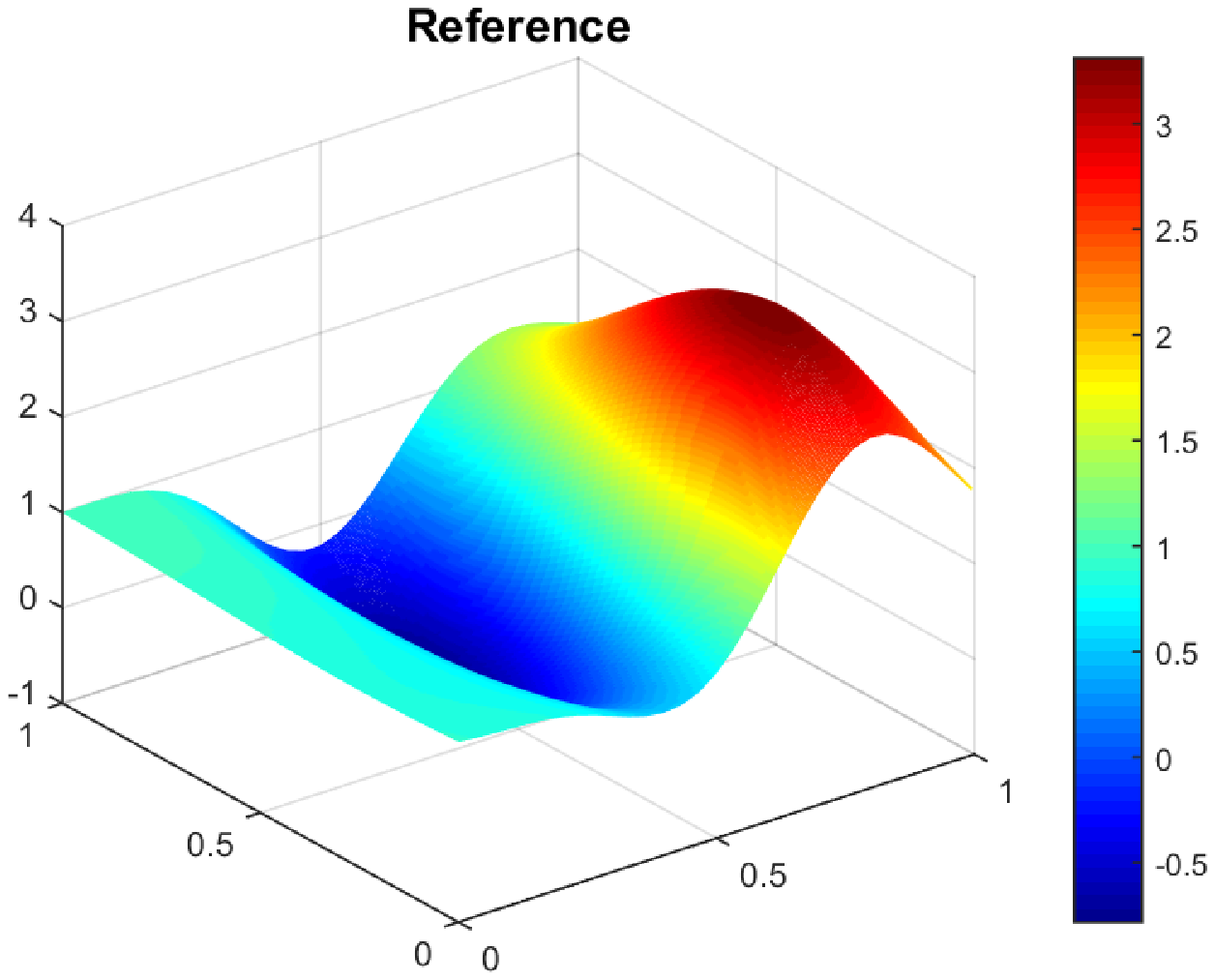}
  \includegraphics[width=1.6in, height=1.7in]{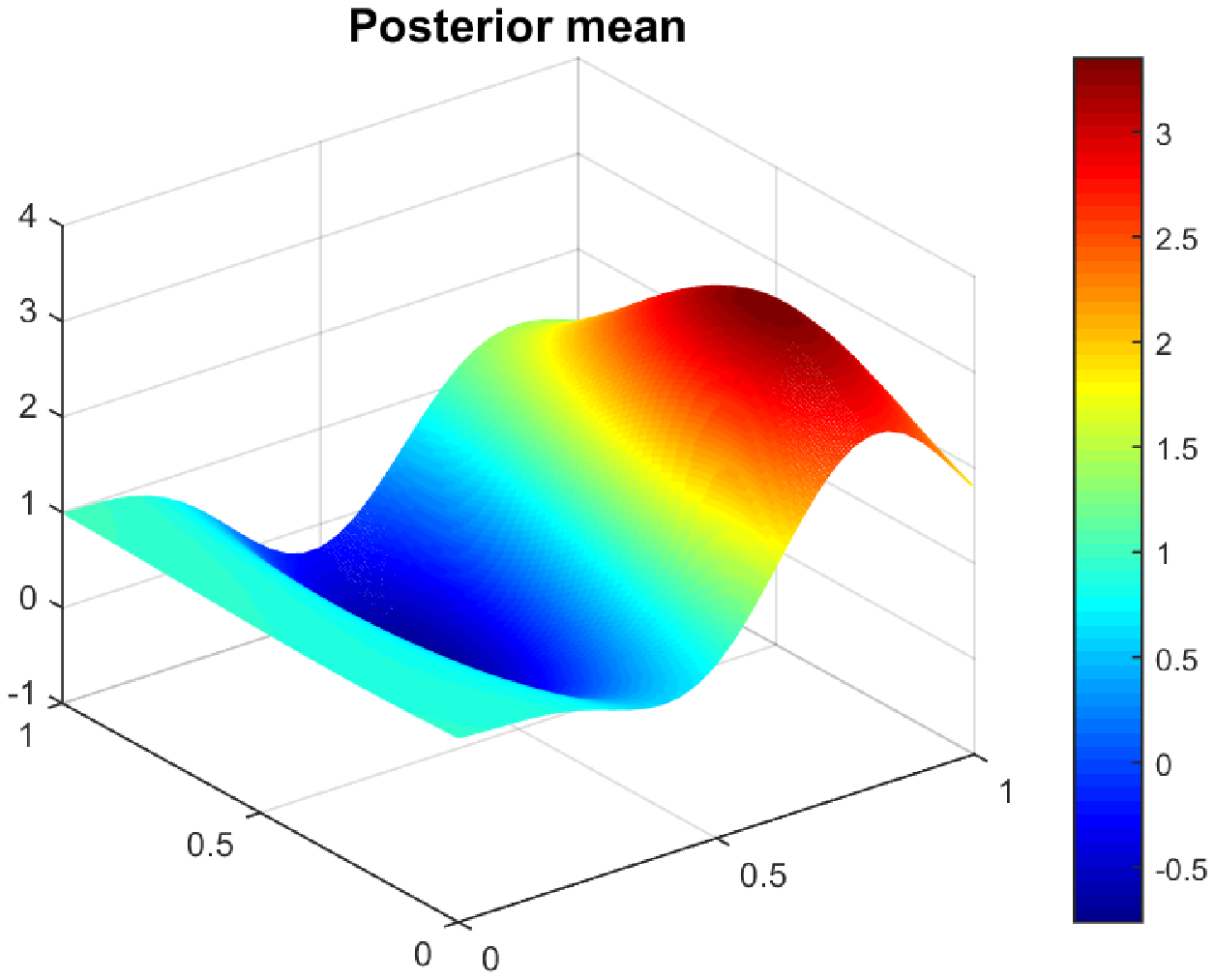}
  \includegraphics[width=1.6in, height=1.7in]{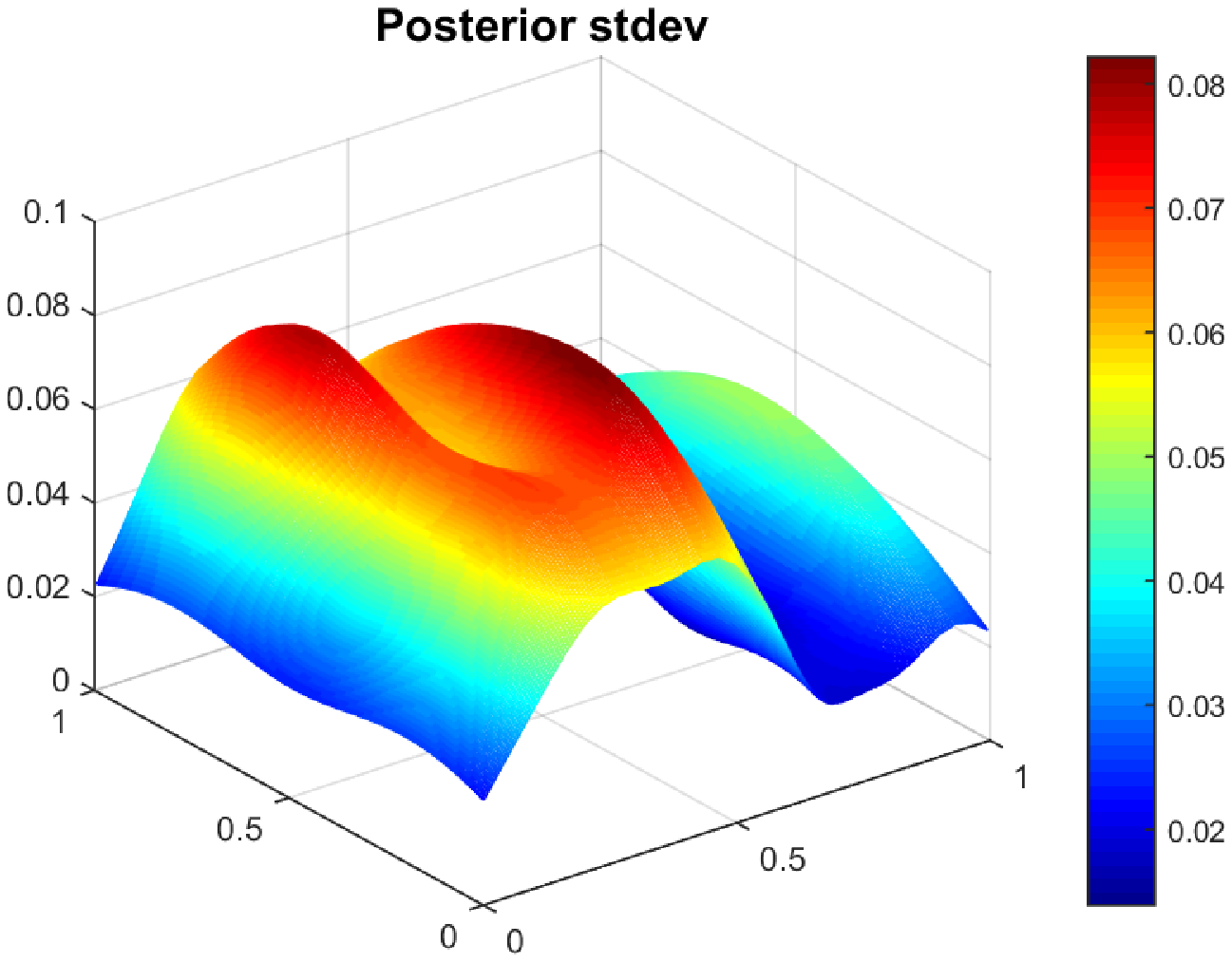}
  \caption{True profile, posterior mean and posterior standard deviation of $\log k(x, \omega)$ (first row) and  $\log q(x, \omega)$ (second row).}
  \label{coef3}
\end{figure}

The forward model is defined  on  $200\times200$ uniform fine  grid. GMsFEM is implemented on  $5\times5$ coarse grid. We set the number of samples in building the snapshot for GMsFEM as $N_{\mu}=15$ and choose $M^i_{\text{snap}}=20$ eigenfunctions in calculating the local snapshot space. GMsFEM with $5$ local multiscale basis functions is used to approximate the forward model in solving the optimization problem (\ref{coarseopt}).  The 2 times standard variance of the approximate sampling distribution for each component of the parameter is plotted in Figure \ref{measure} (right). We adjust $\nu$ to be 1 for the first 10 components and set $\nu$ to be the $2\sigma_\text{OLS}\sqrt\delta$  for the last 5 components.

The gPC order is set as $N=3$ for LS-SCM.  We run 5 Markov chains using intermediate-based surrogate.   The result is shown in  Figure \ref{coef3}. For the logarithmic specific field, the standard variation is small in the  physical domain, and uncertainty lies among the low-lying domain of the field and the junction of the peak and valley of the logarithmic  specific field. The standard variation is symmetrically distributed for the logarithmic permeability field, as the reference field appears to be fluctuate and its uncertainty concentrates at peaks and valleys.

\begin{figure}
\centering
 \includegraphics[width=2in, height=1.7in]{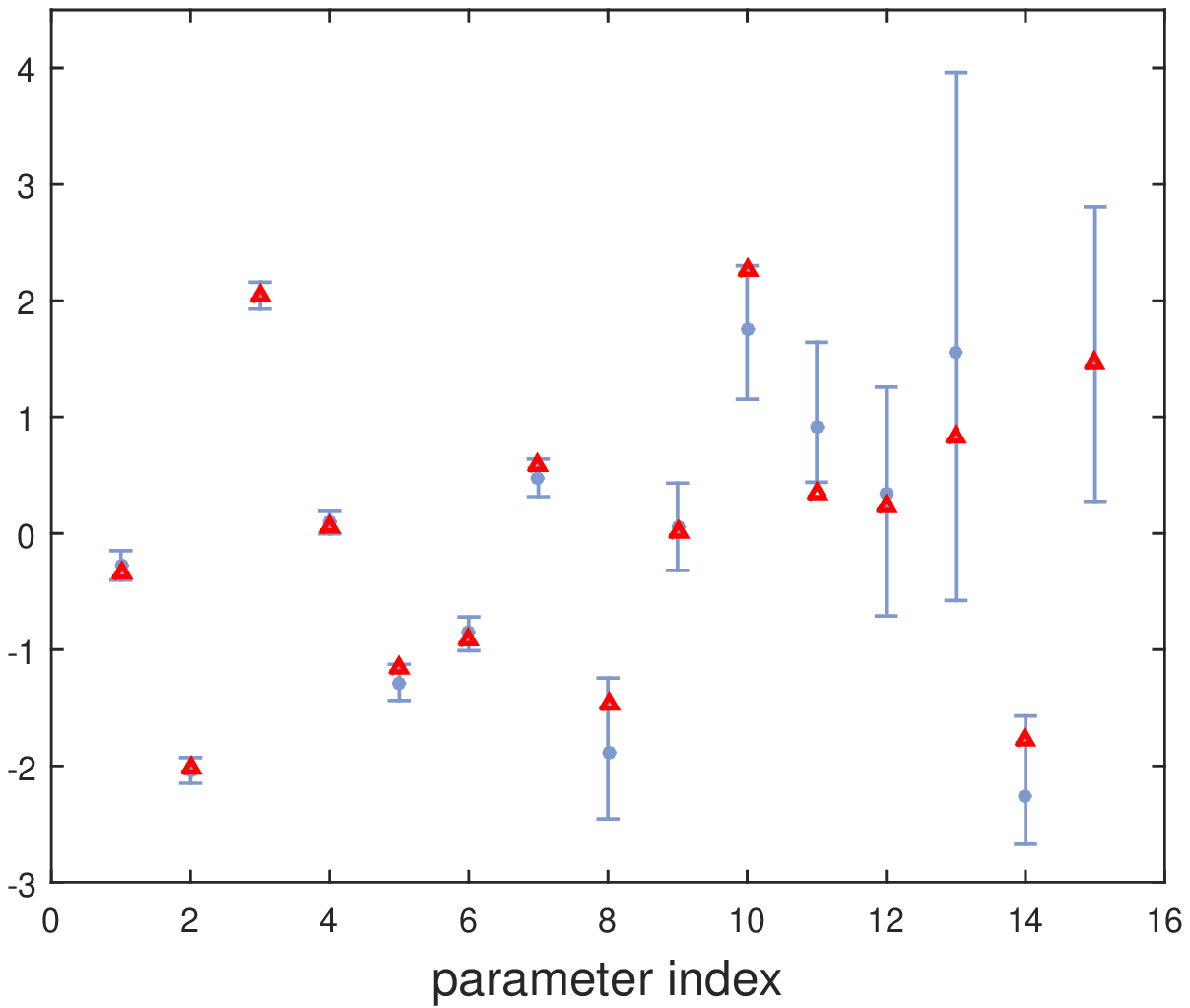}
  \includegraphics[width=2in, height=1.7in]{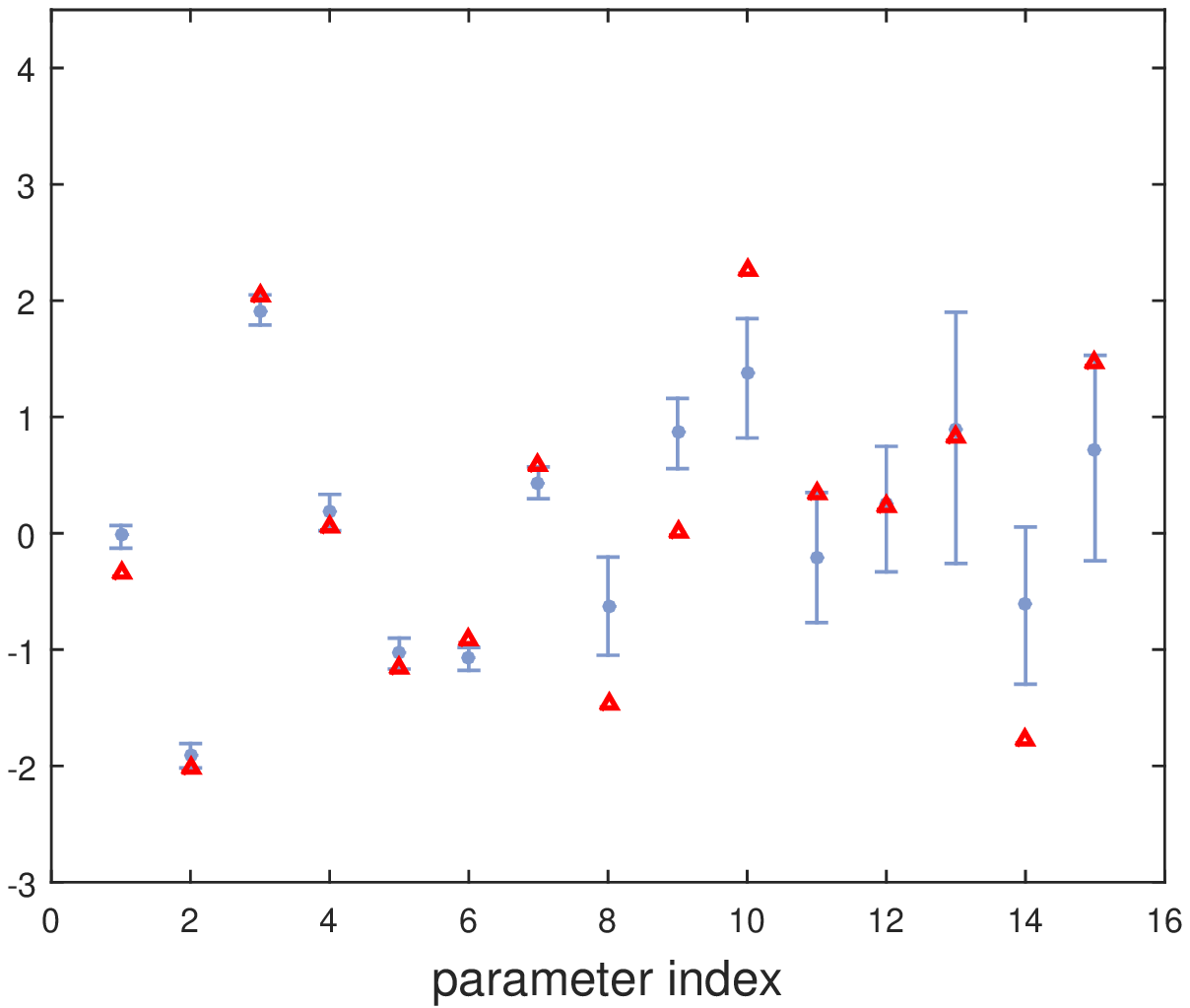}
\caption{Point estimates (gray filled circle), reference (red unfilled triangle) and 95\% credible intervals for $z$ , (left) induced by intermediate-based surrogate, and (right) prior-based surrogate}\label{credi}
\end{figure}

The 95\% credible intervals for parameter $z$ by posterior samplers with respect to the intermediate-based surrogate and the prior-based surrogate  are plotted in Figure \ref{credi}. We note that both intervals are tight for the first 7 parameters and the uncertainty concentrates on the parameters associated with the permeability field for the two types of surrogate. This is due to the
 small correlation length $l_{y, Y_1}=0.01$ and the exponential covariance kernel in generating the reference field for the permeability field. For the intermediate-based surrogate, the credible values are symmetric except for $z_8$ and $z_{14}$, and almost all parameters are comprised within the 95\% uncertainty range, while for the prior-based surrogate, there are 5 out of 15 reference points lie outside of the 95\% credible interval.  This  implies that the intermediate-based surrogate performs better than the prior-based surrogate.

\begin{figure}
\centering
\includegraphics[width=2in, height=1.7in]{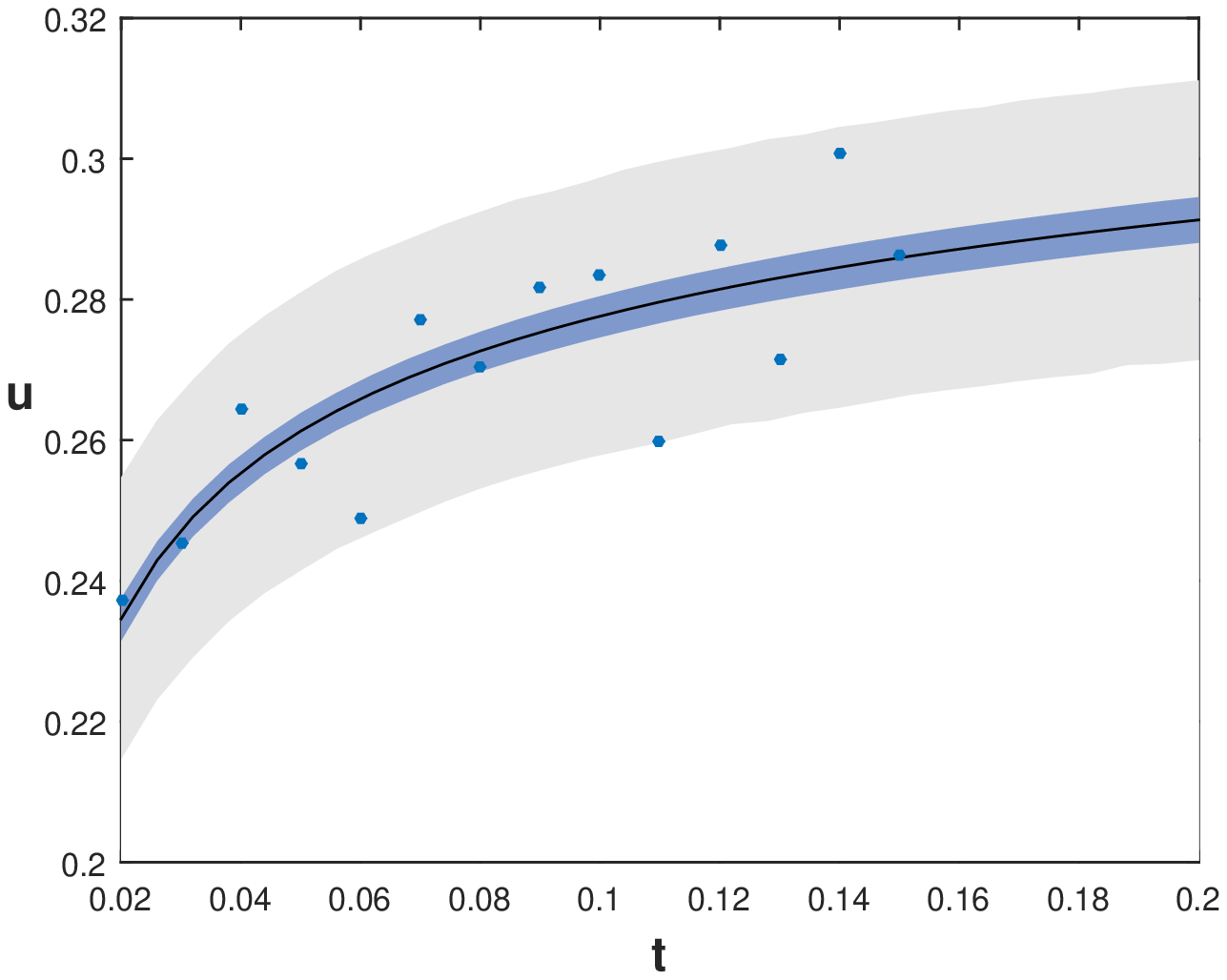}
 \includegraphics[width=2in, height=1.7in]{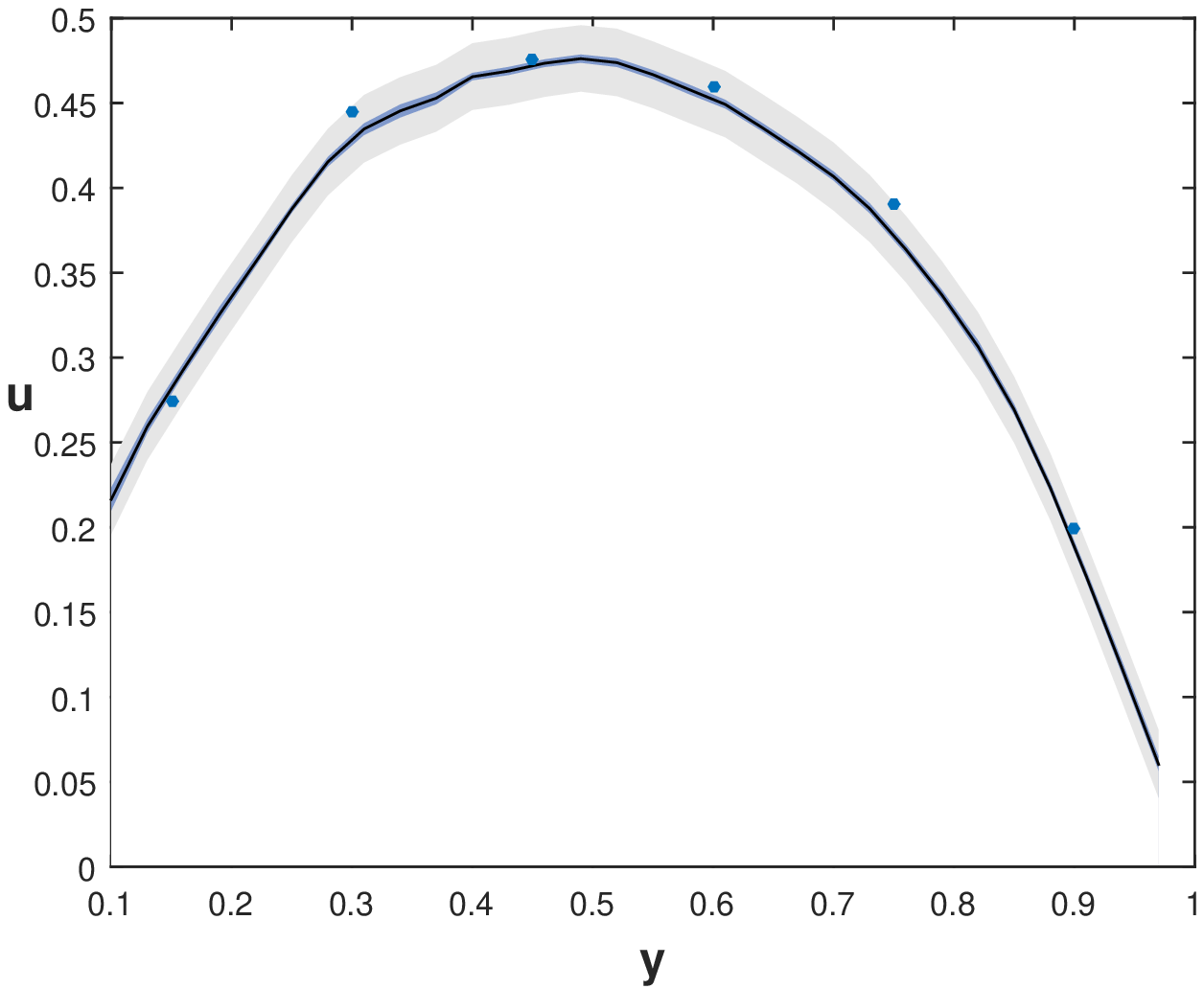}
\caption{Data, point estimates, and 95\% credible and prediction intervals produced by the Bayesian analysis for $u((0.15,0.3);t)$ (left), and $u((0.6,y);0.1)$ (right).}
\label{pret}
\end{figure}

We also construct the 95\% credible intervals for model response at $u\big((0.15, 0.3); t\big)$ and $u\big((0.6, y); 0.1\big)$, with 25000 posterior samplers and realizations constructed by the intermediate-based surrogate. This uncertainty is added to the estimated error variance to construct prediction intervals. As illustrated in Figure \ref{pret}, measurement data are almost contained in the predictive intervals. For $u\big((0.15, 0.3); t\big)$, both the credible interval and prediction interval get wider as time moves on.  This means that as the uncertainty from input $z$ propagates  with respect to  time, the uncertainty associated with the model fit and predictions grows. The uncertainty mainly concentrates on $y\in (0.3, 0.7)$ for $u((0.6, y); 0.1)$, and as $y$ gets closer to the end points, both credible interval and prediction interval  becomes tight, this is due to the deterministic Dirichlet boundary condition  at $y=0$ and $y=1$.

\subsection{Identification of a mixed random field}
\label{sec4-3}
In this subsection,  we consider the fractional diffusion equation with the same boundary condition  as in Subsection  \ref{rec-k},    but the unknown coefficient $k(x,\omega)$ has the representation
  \[
  k(x,\omega)=\sum_i^m k_i(x,\omega)\bff{I}_{\Omega_i}(x),
  \]
where $\bff{I}_{\Omega_i}$ is an indicator function on $\Omega_i\subset\Omega$, and each $k_i(x,\omega)$ has the truncated KLE form.
We assume $m=2$ and the indicator functions are known.  It is assumed that there is a high permeability layer in the middle and low permeability layer in the two ends, the high permeable layer has correlation lengths $l_{x,Y_1}=0.1$, $l_{y,Y_1}=0.4$ and $\sigma_{Y_1}=1$, $\bb{E}[Y_1(x,\omega)]=0$, the low permeable layer has correlation lengths $l_{x,Y_2}=0.2$, $l_{y,Y_2}=0.2$ and $\sigma_{Y_2}=0.4$, $\bb{E}[Y_2(x,\omega)]=0$. We assume they are correlative with the correlation coefficient $\rho=0.8$, i.e., the truncated random field has the form
  \begin{eqnarray*}
  Y_1(x,\omega) &=& \log{k_1(x,\omega)}=\sigma_{Y_1}\sum_{i=1}^{N_1} \sqrt{\lambda^1_i}e^1_i(x)z^1_i(\omega), \\
  Y_2(x,\omega) &=& \log{k_2(x,\omega)}=\sigma_{Y_2}\bigg(\rho\sum_{i=1}^{N_1} \sqrt{\lambda^1_i}e^1_i(x)z^1_i(\omega)+
  \sqrt{1-\rho^2}\sum_{i=1}^{N_2} \sqrt{\lambda^2_i}e^2_i(x)z^2_i(\omega)\bigg),
  \end{eqnarray*}
where $Y_1$ is the low permeability field and $Y_2$ is the high permeability field, both eigenfunctions $\{e^1_i(x)\}_i^{N_1}$ and $\{e^2_i(x)\}_i^{N_2}$ and are generated from  Gaussian type kernel. The truncation term in the KLE is set as $N_1=9, N_2=15$ for $Y_1$ and $Y_2$, respectively. Measurements are taken at time  instances $0.02:0.01:0.21$ and the measured locations are uniformly distributed on $\Omega$ with stepsize 35/200.

\begin{figure}
  \centering
  \includegraphics[width=1.6in, height=1.7in]{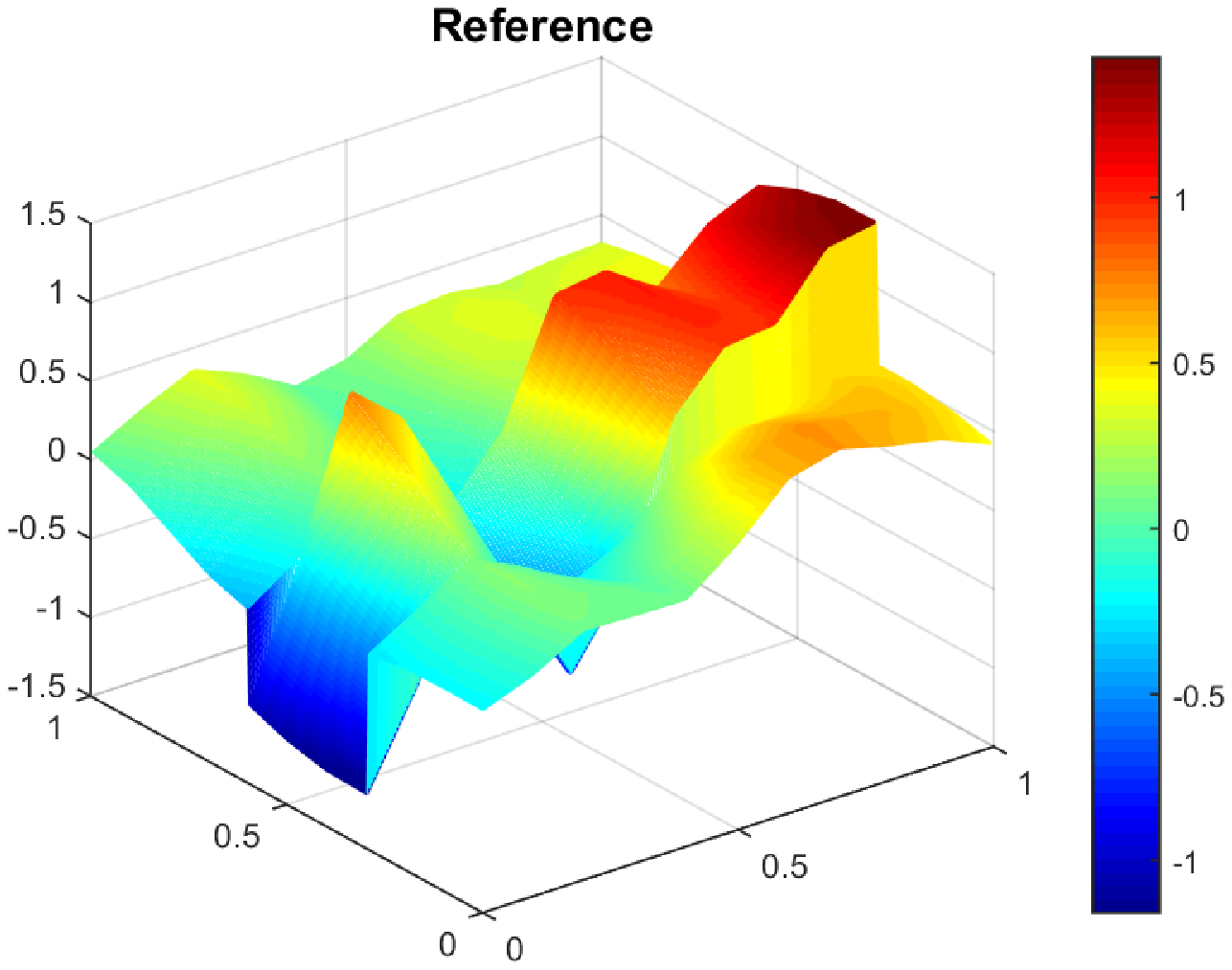}
  \includegraphics[width=1.6in, height=1.7in]{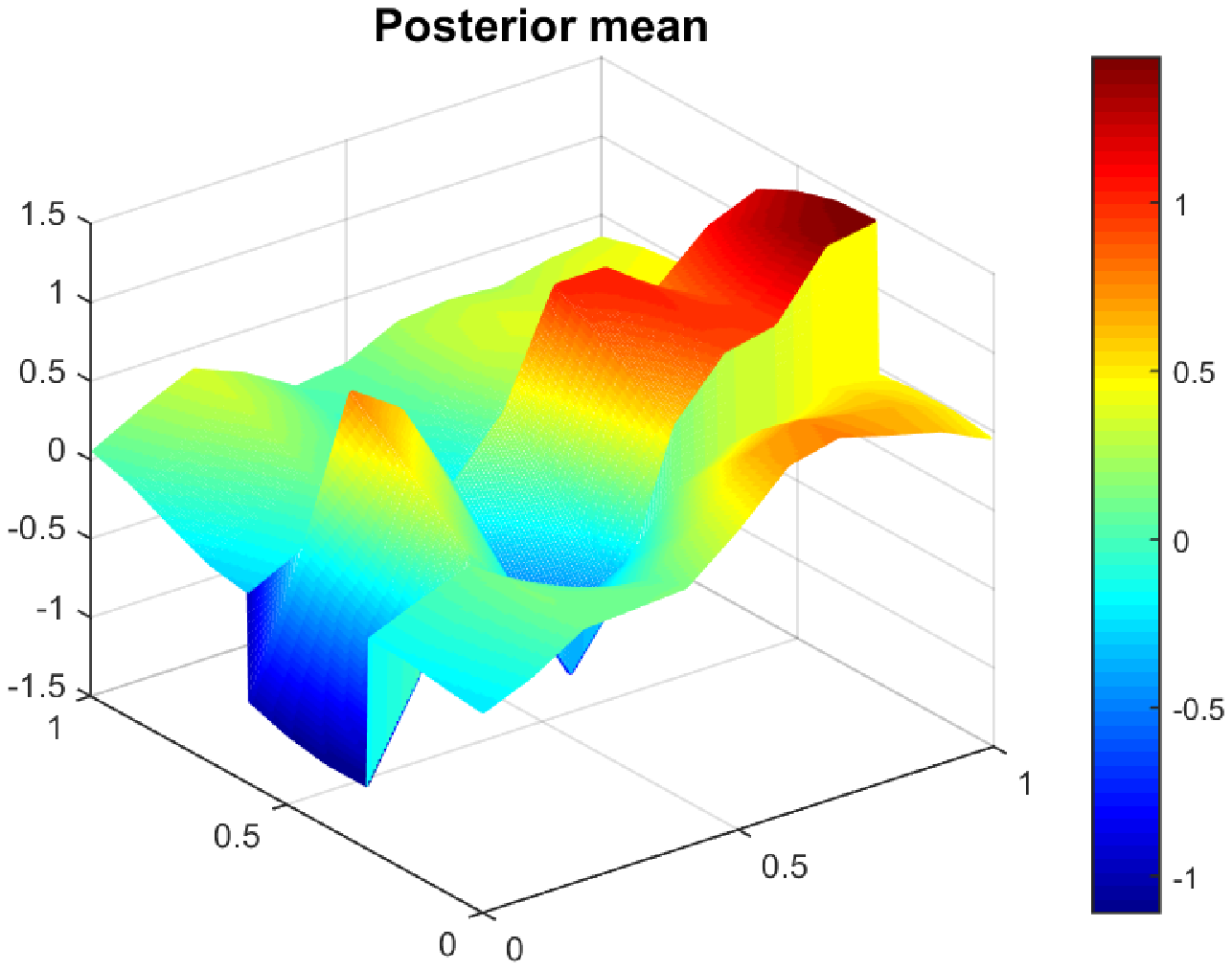}
  \includegraphics[width=1.6in, height=1.7in]{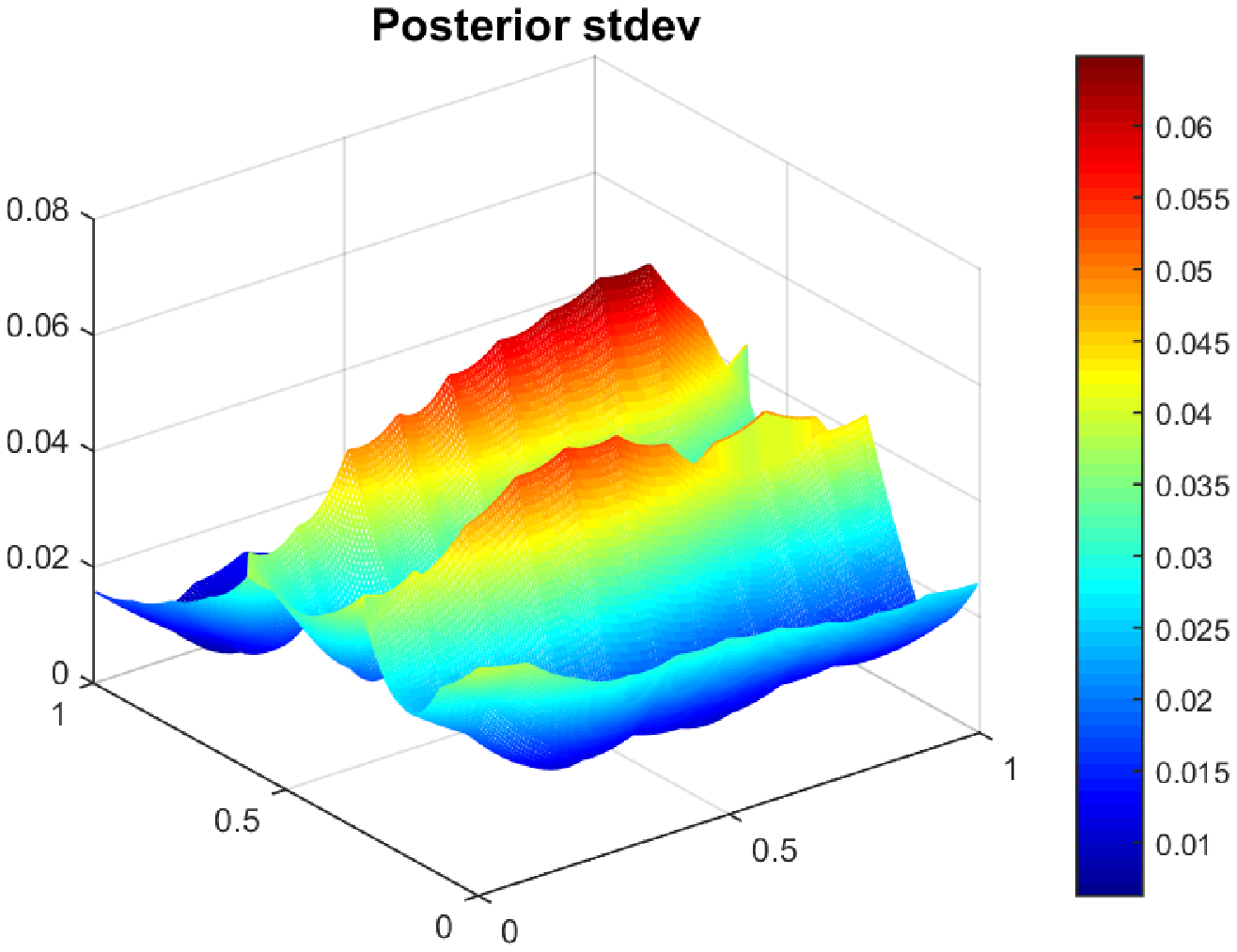}
  \caption{True profile, posterior mean, posterior standard deviation of  $\log k(x, \omega)$.}\label{coef4}
\end{figure}

The forward model is defined  on  $200\times200$ uniform grid, and GMsFEM  is implemented on $5\times 5$ coarse grid.   GMsFEM with $7$ local multiscale basis functions (i.e., $M_c=7$)
is used to solve   optimization problem (\ref{coarseopt}).   As we have known some degree of accuracy can be achieved with intermediate-based surrogate even when gPC order is low, we set the order of gPC as $N=2$ in LS-SCM  surrogate. Thus we need to  solve the forward GMsFE  model for 650 times. Because the bounds $2\sigma_\text{OLS}\sqrt\delta$ are smaller than 1 for all parameter components,  we construct the gPC surrogate with respect to intermediate distribution $U(z_\text{OLS}-1, z_\text{OLS}+1)$. The logarithmic true profile, posterior mean and posterior standard variation are illustrated in Figure \ref{coef4}, which shows that the uncertainty mainly lies on junction of high and low permeable layer.

\begin{figure}
\centering
\includegraphics[width=2in, height=1.7in]{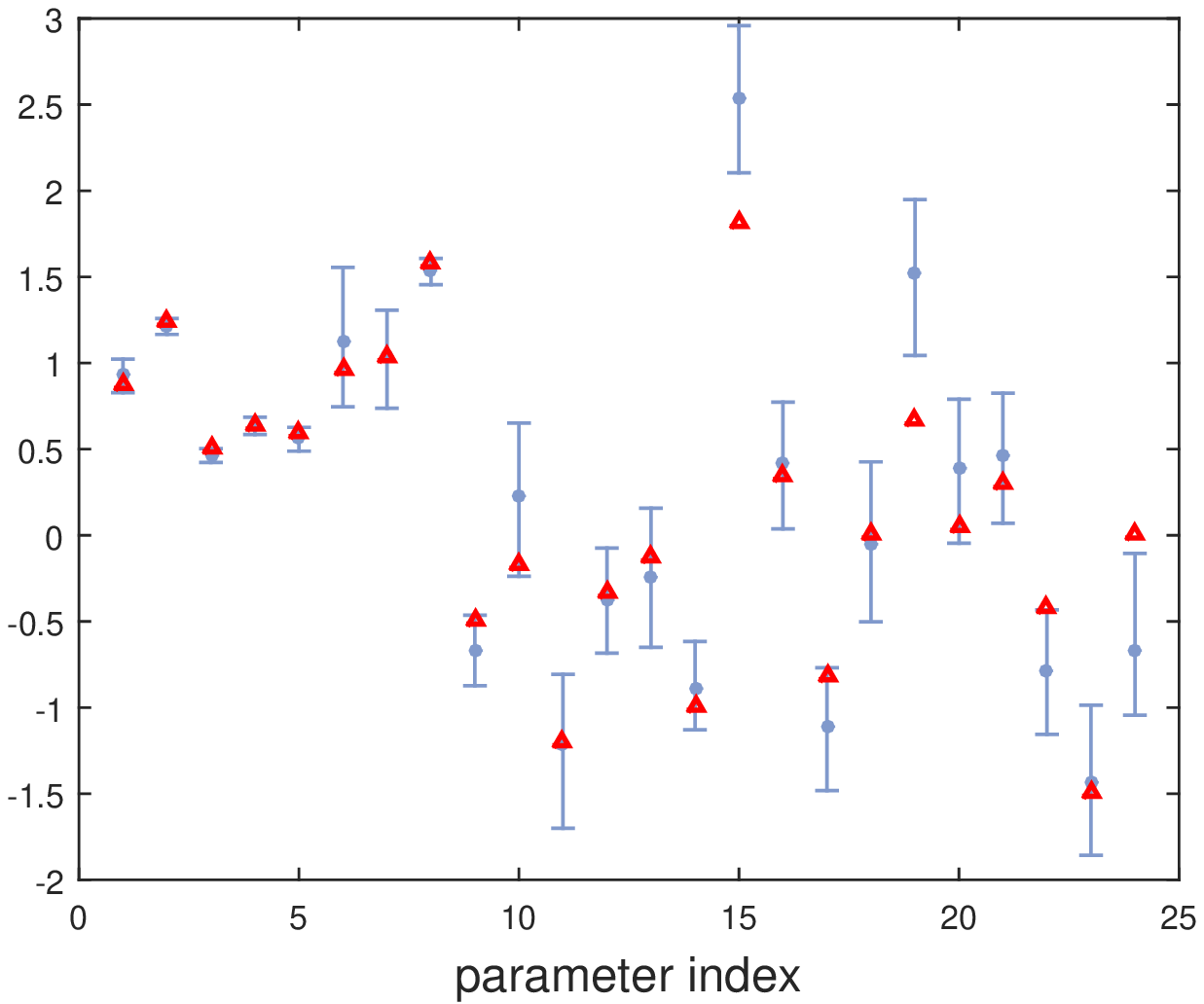}
 \includegraphics[width=2in, height=1.7in]{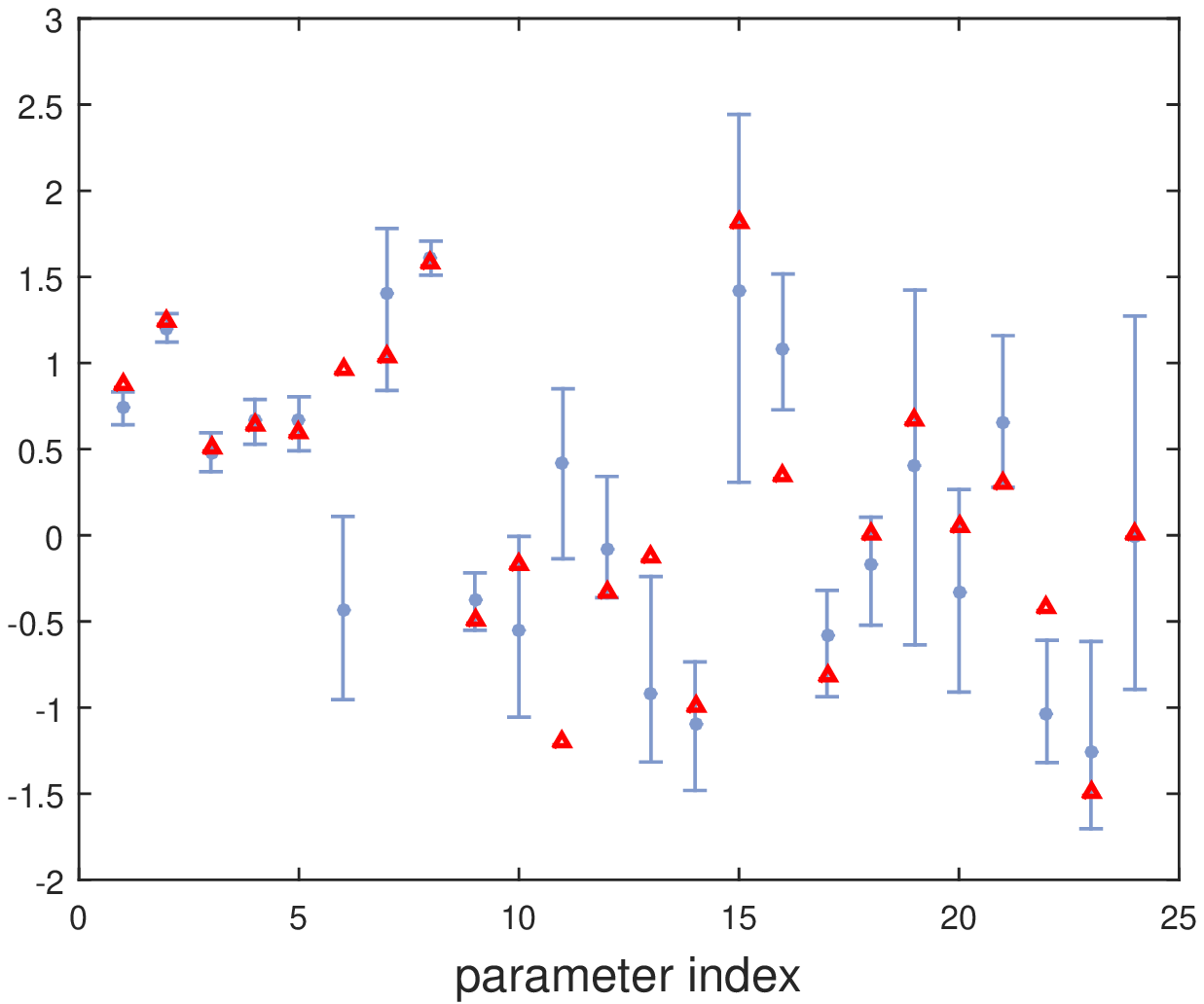}
 \caption{Point estimates (gray filled circle), reference (red unfilled triangle) and 95\% credible intervals for $z$ , the left one is induced by intermediate-based surrogate, and the right one is induced by prior-based surrogate.}
 \label{credi2}
\end{figure}

Figure $\ref{credi2}$ depicts the 95\% uncertainty bands of parameter $z$ obtained by MCMC samplers. There are 5 out of 25 reference points lie outside of the 95\% credible interval with respect to prior-based surrogate while the proportion of points lying outside of the credible interval is 3 out of 25 with respect to intermediate-based surrogate. Though the difference in proportion is small, we can see the magnitude of deviation is much different, especially for $z_6$ and $z_{11}$ as shown in Figure \ref{credi2}. We emphasize that the intermediate distribution  gives  us  the intervals where posterior has high possibility, most of the unimportant  support of posterior has been excluded out.  Thus  the intermediate-based surrogate performs better than the prior-based one when gPC order is low.

\begin{figure}
\centering
\includegraphics[width=1.6in, height=1.7in]{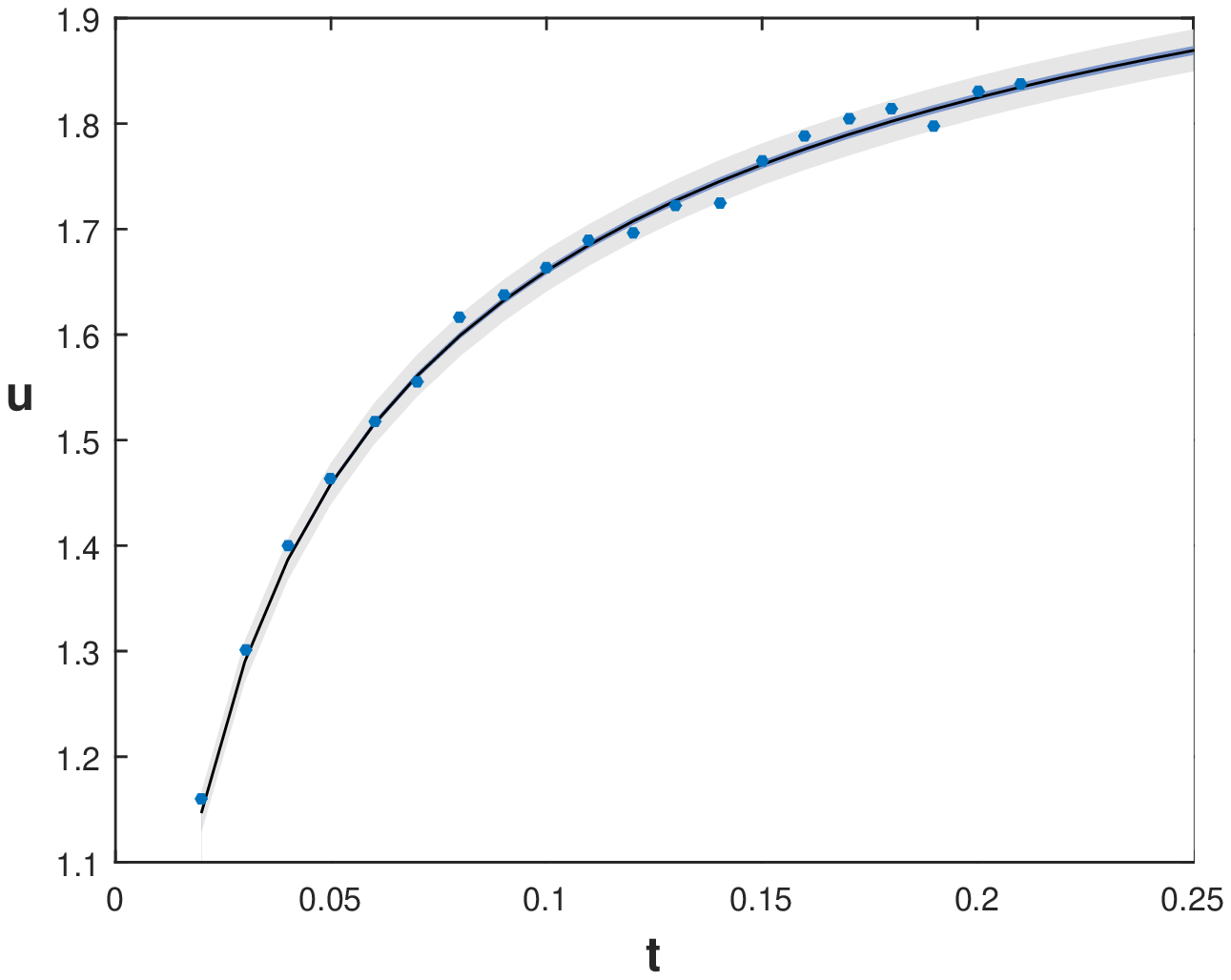}
 \includegraphics[width=1.6in, height=1.7in]{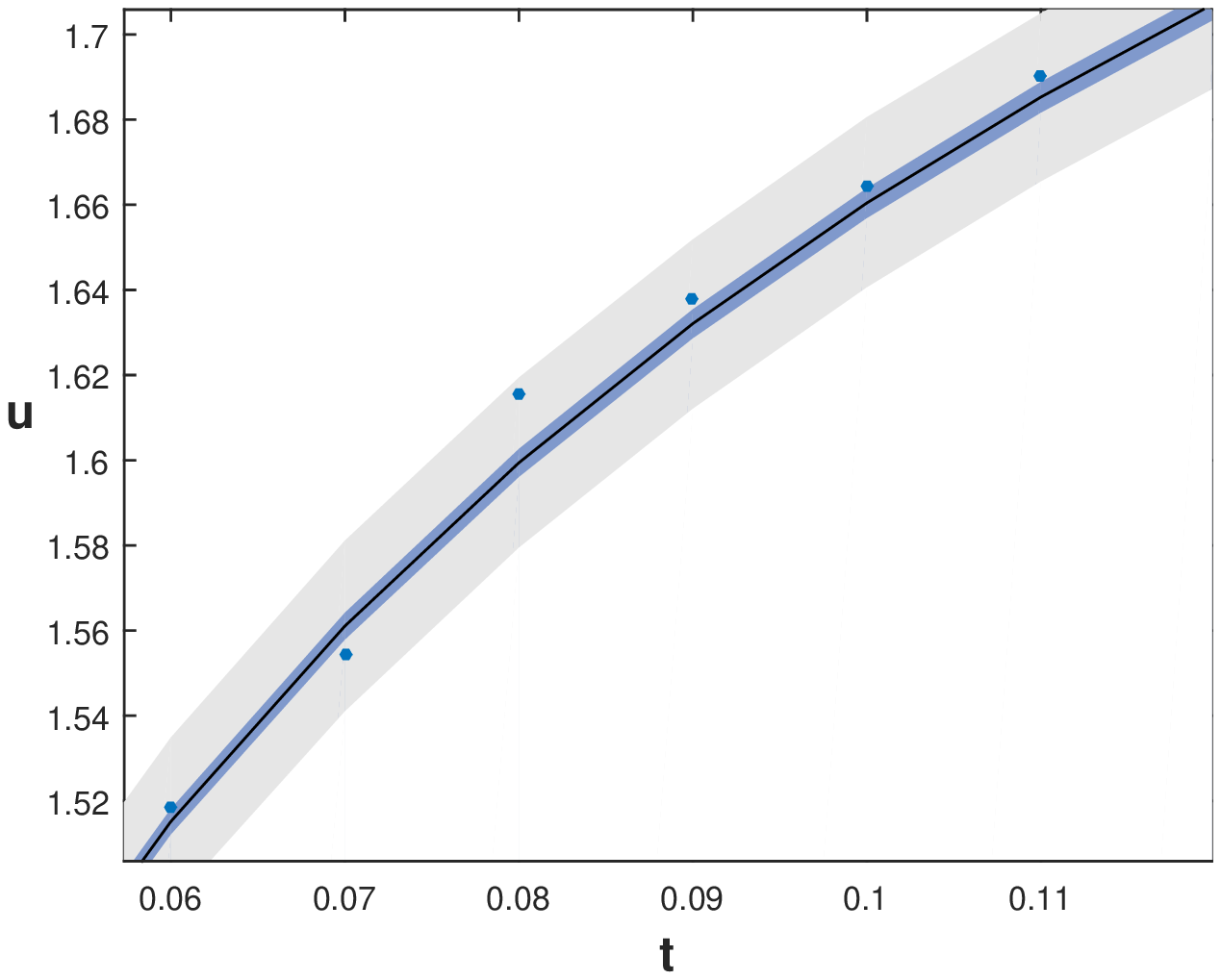}
 \includegraphics[width=1.6in, height=1.7in]{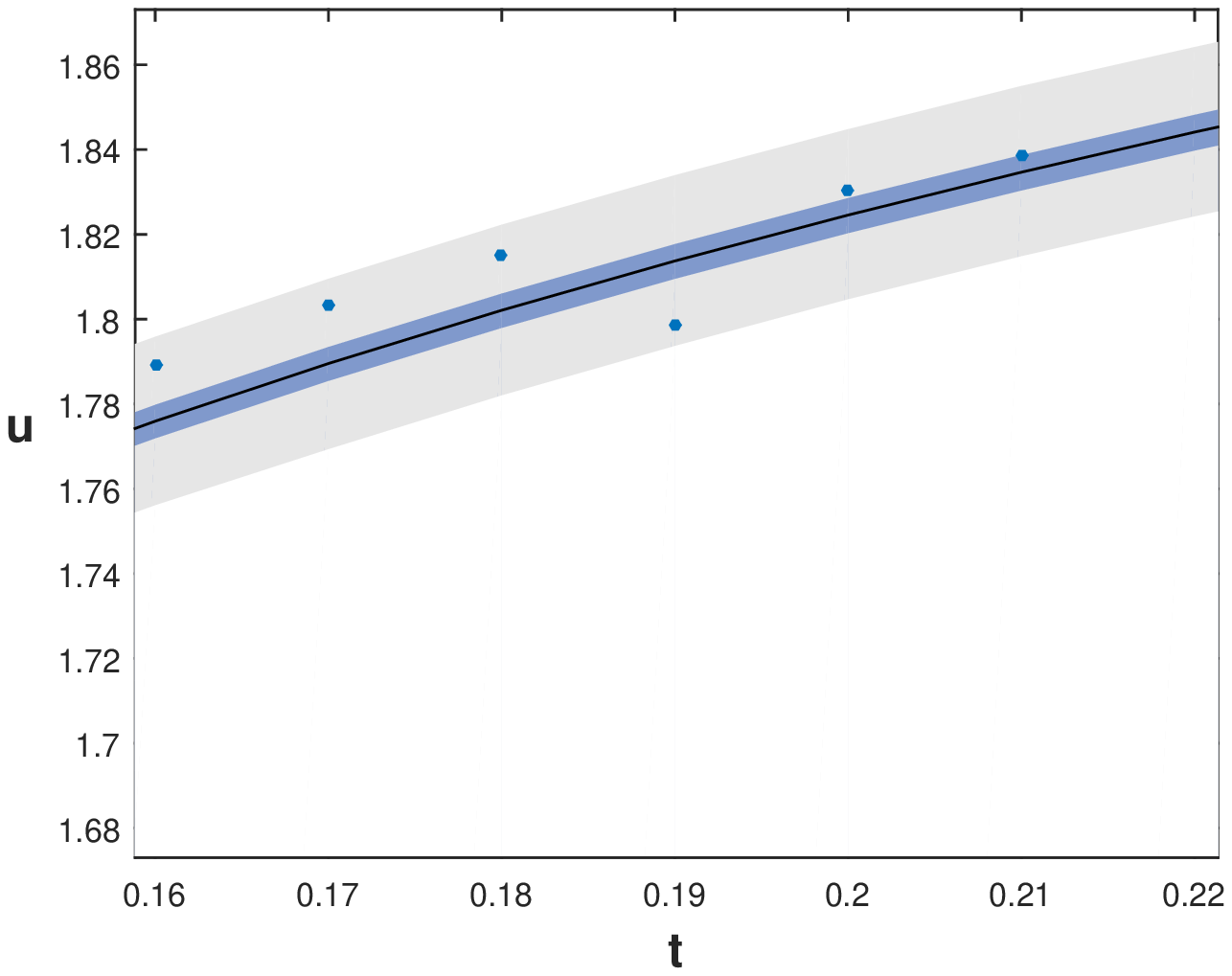}
 \caption{Data, point estimates, and 95\% credible and prediction intervals for $u((0,0);t)$ produced by the Bayesian analysis: (left) full time interval and (middle and right) partial time interval to illustrate different interval width for the beginning and end time.}
 \label{pre-t}
\end{figure}

\begin{figure}
\centering
\includegraphics[width=1.6in, height=1.7in]{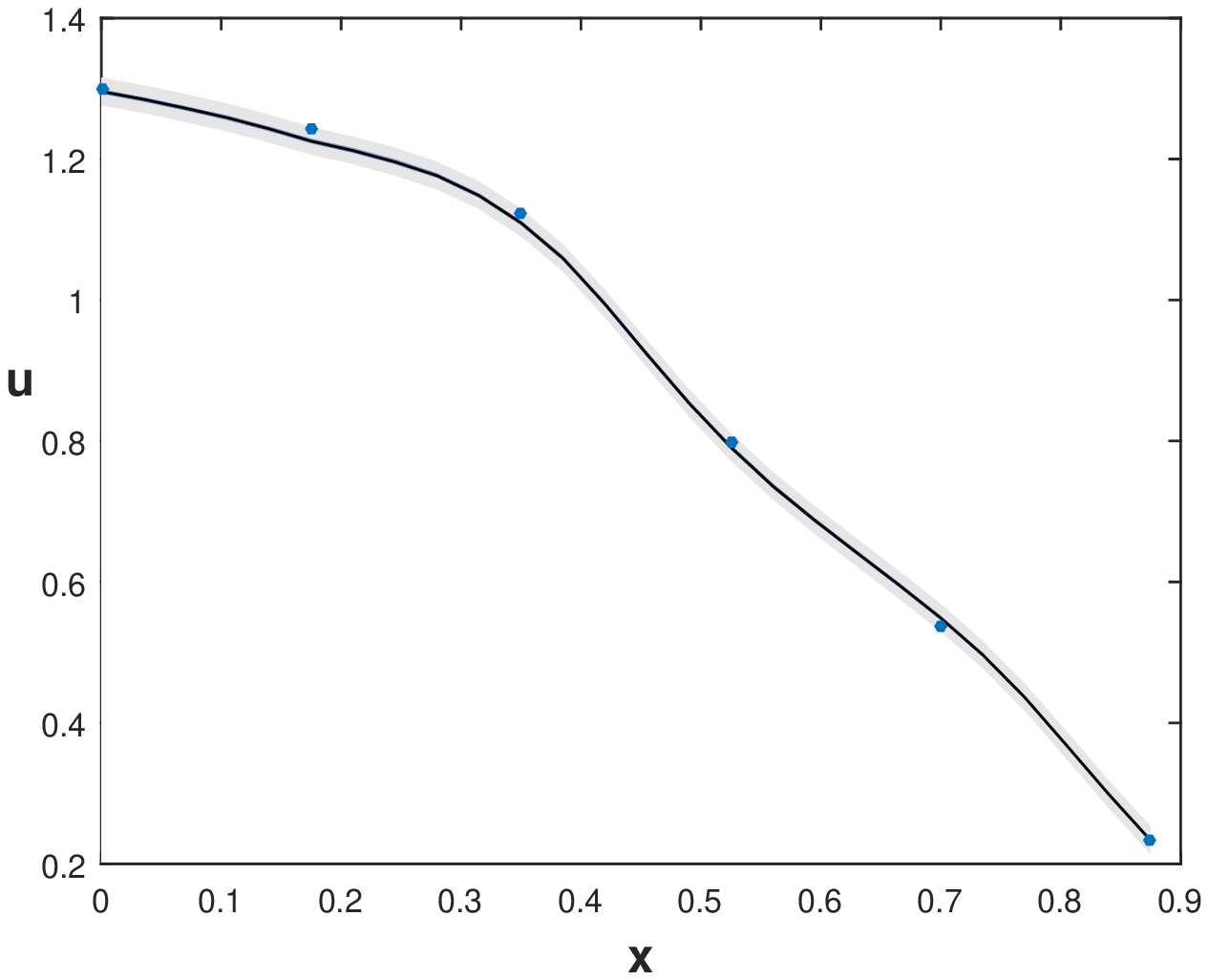}
 \includegraphics[width=1.6in, height=1.7in]{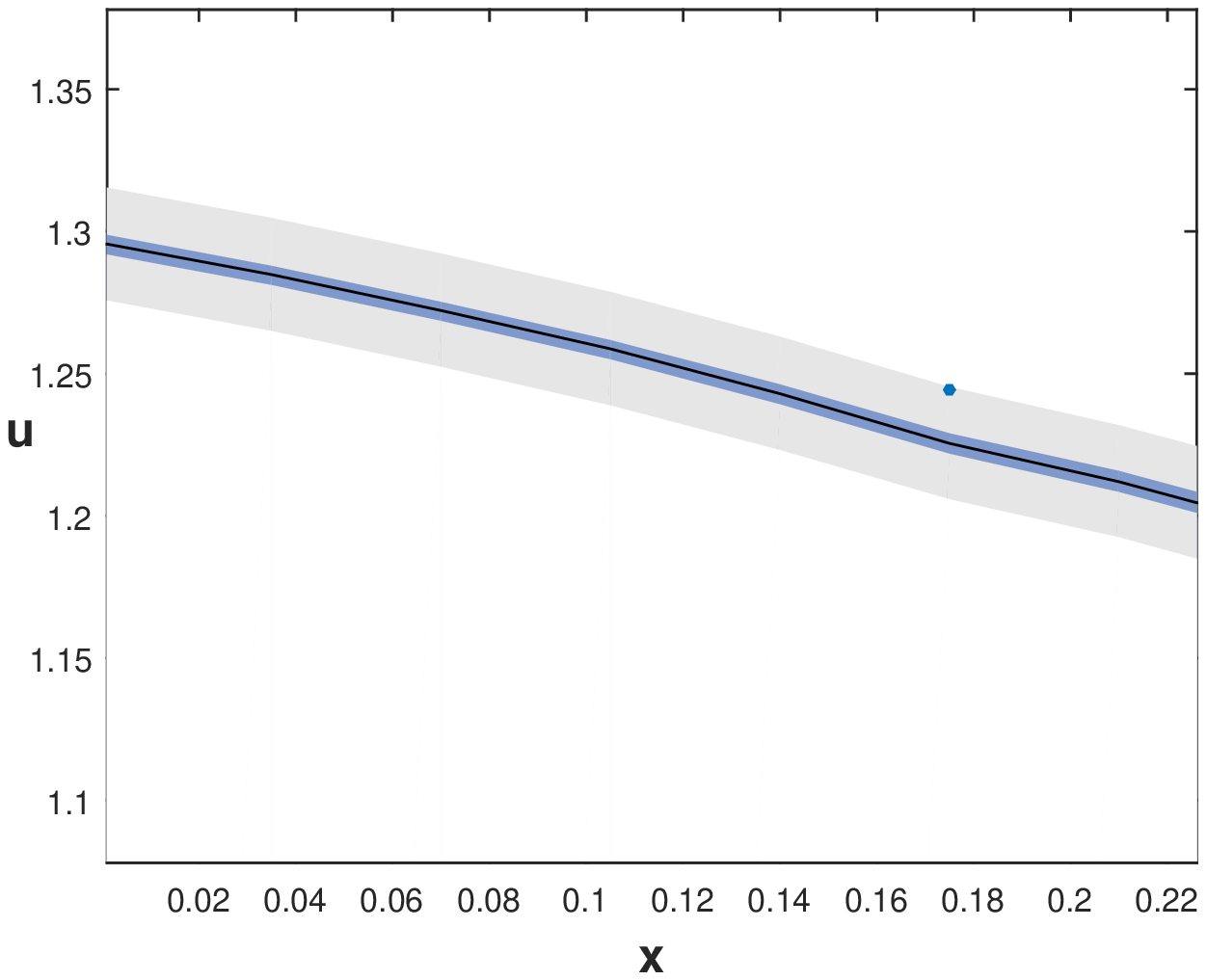}
 \includegraphics[width=1.6in, height=1.7in]{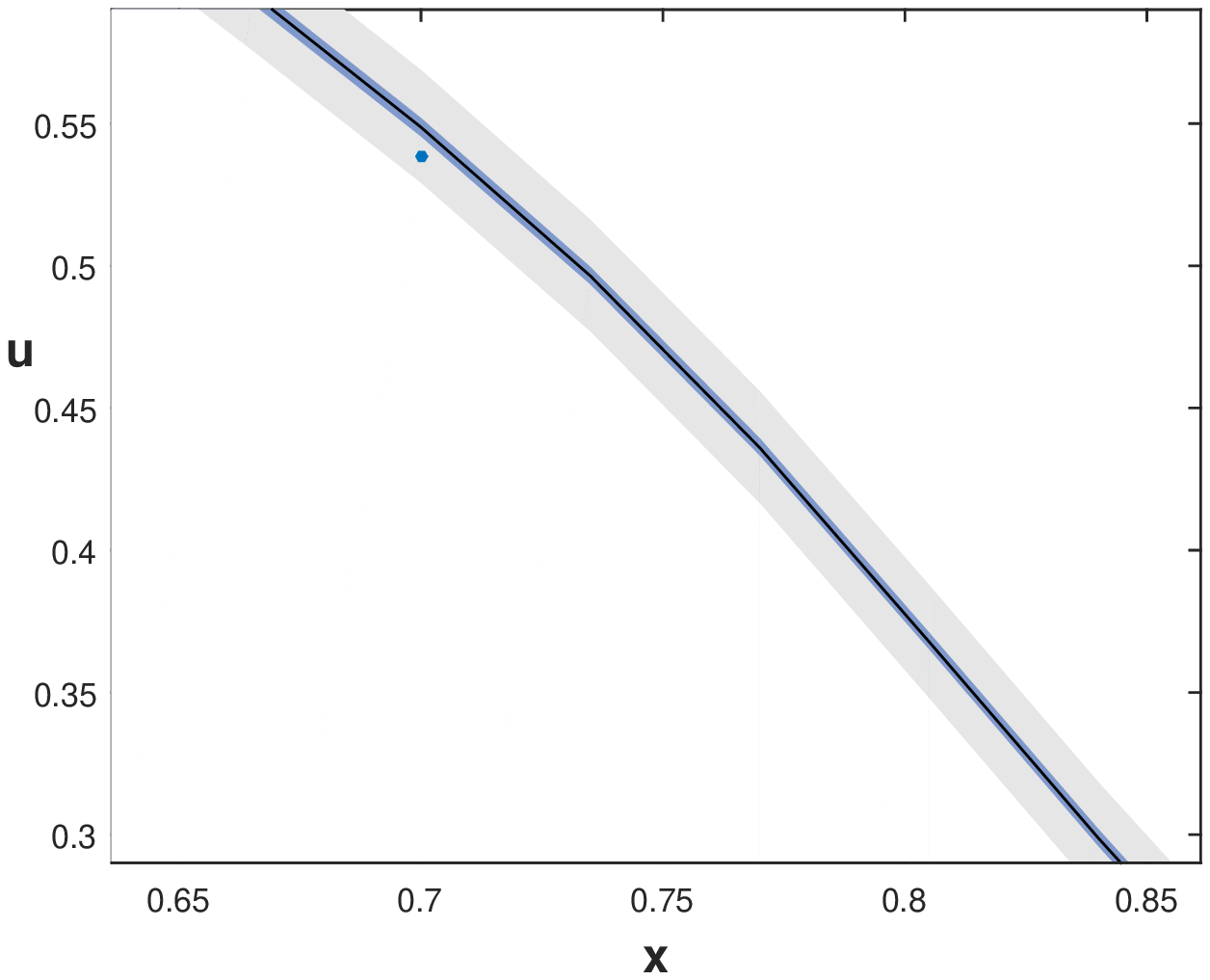}
 \caption{Data, point estimates, and 95\% credible and prediction intervals for $u((x,0.7);0.1)$ produced by the Bayesian analysis, the solid line is the posterior mean: (left) full physical $y=0.7$ interval and (middle and right) partial  physical interval to illustrate different interval width for different $x$.}
 \label{pre-x}
\end{figure}

To construct 95\% credible intervals for model response, we use samplers from the parameter chains to construct 25000 realizations by the intermediate-based surrogate. The prediction intervals are constructed by accounting for measurement errors.  The credible  interval and prediction interval, along with the point estimates, and synthetic data, are illustrated for $u\big((0,0); t\big)$ and $u\big((x, 105/200); 0.1\big)$ in Figure \ref{pre-t} and Figure \ref{pre-x}, respectively. We note that both intervals are symmetric with respect to the point estimate, and both credible intervals are tight for $u\big((0,0); t\big)$ and $u\big((x, 105/200); 0.1\big)$. For $u\big((0,0); t\big)$, we can observe from Figure \ref{pre-t} (middle) (right) that the uncertainty associated with both the model fit and predictions grows with respect to  time, while for $u\big((x, 105/200); 0.1\big)$, the corresponding uncertainty decreases as $x$ gets closer to $x=1$, which is the deterministic  Dirichlet boundary condition.

\section{Conclusion}

This paper presented a two stage strategy on constructing surrogate model by GMsFEM and LS-SCM. In the first stage, we used a reduced order model to approximate the forward model, and the corresponding misfit-to-observed nonlinear least square problem is solved by using the regularizing Levenberg-Marquart algorithm and GMsFEM.
 We then constructed the intermediate distribution based on the approximate sampling distribution, and used  GMsFEM and LS-SCM to construct the representation for the model response based on the intermediate distribution.   The DREAM$_\text{ZS}$ algorithm was  employed to explore the surrogate posterior density, which is incorporated by the surrogate likelihood and the intermediate distribution. It showed  that the proposed method leads to the approximate posterior as accurate as the one derived directly based on the original prior, with lower gpc order.

We have combined the deterministic and statistical methods together to solve inverse problems.  We obtained not only the point estimate or confidence interval but also the statistical properties of the unknowns. We solved   optimization problem by gradient based method for a coarse model, and constructed the intermediate distribution artificially, which do not need to approximate the posterior accurately  but contains the significant  region or support of the posterior.  The intermediate distribution is incorporated with the likelihood to be explored, and the acceptance rate is
considerably improved as well.


\end{document}